\documentclass[12pt]{amsart}
\usepackage[margin=1in]{geometry}
\usepackage[utf8]{inputenc}
\usepackage{amsmath}
\usepackage{amssymb}

\usepackage[pdftitle={DisquesPFH4},
  pdfauthor={Samuel Bronstein},
  pdffitwindow=true,
  breaklinks=true,
  colorlinks=true,
  urlcolor=blue,
  linkcolor=red,
  citecolor=red,
  anchorcolor=red]{hyperref}
\usepackage{titleps}
\usepackage{xcolor}
\definecolor{links}{RGB}{50,0,200}
\definecolor{hyperrefcolor}{rgb}{50,0,250}
\hypersetup{citecolor=blue,linkcolor=links}

\usepackage{subfiles}
\usepackage{pgfplots}
\pgfplotsset{compat=1.15}
\usepackage{mathrsfs}

\numberwithin{equation}{section}
\numberwithin{table}{section}
\numberwithin{figure}{section}
\newcommand{\B}{\mathbb B}

\newcommand{\E}{\mathbb E}
\newcommand{\N}{\mathbb N}
\newcommand{\R}{\mathbb R}

\renewcommand{\H}{\mathbb H}
\renewcommand{\S}{\mathbb S}

\newcommand{\calH}{\mathcal H}

\newcommand{\calP}{\mathcal P}

\newcommand{\calV}{\mathcal V}

\newcommand{\eps}{\varepsilon}

\newcommand{\Vol}{\text{Vol}}

\newcommand{\one}{\mathbf{1}}

\theoremstyle{plain}
\newtheorem{theorem}{Theorem}[section]
\newtheorem*{theorem*}{Theorem}
\newtheorem{lemma}[theorem]{Lemma}
\newtheorem{proposition}[theorem]{Proposition}

\newtheorem{corollary}[theorem]{Corollary}
\newtheorem*{corollary*}{Corollary}

\theoremstyle{definition}
\newtheorem{definition}[theorem]{Definition}

\newtheorem{example}[theorem]{Example}
\newtheorem{notation}[theorem]{Notation}

\theoremstyle{remark}
\newtheorem*{remark}{Remark}

\title{Dominating metric measure spaces by radial spaces}
\author{samuel bronstein}
\date{2025}

\begin{document}
\maketitle
\begin{abstract}
	This paper is devoted to a kind of rearrangement of functions on $\rm{CD}(k,n)$-spaces,
	which satisfy a Polyà--Szegö type inequality.
	We use this rearrangement to prove the validity of a Moser--Trudinger type inequality on
	a wide class of metric measure spaces satisfying a $\rm{CD}(k,n)$-curvature-dimension inequality.
	As a consequence, we give a characterization, among manifolds with
	lower bounded Ricci curvature, of those admitting a Moser--Trudinger type inequality.
\end{abstract}
\tableofcontents
\section{Introduction}
In this paper we consider the Moser--Trudinger inequality on metric measure spaces satisfying a $\rm{CD}(k,n)$-property.
a class of metric measure spaces with lower bounded synthetic Ricci curvature.

First proven on the sphere~$\S^2$ by Trudinger in \cite{Tru67}, the Moser--Trudinger inequality
has since been extended by Moser~\cite{Mos71} to n-dimensional spheres,
and more recently to a wide variety of compact and noncompact manifolds see \cite{Li01,Li05,MS10,Yan12,MST13,Kri19,LL21}.
It always amounts to, under some conditions on the $L^n$-norm of the gradient
of a function, to get a bound on the integral of $e^u-P_n(u)$,
where $P_n(t)=\sum_{j=0}^{n-2}\frac{t^j}{j!}$ is the first n terms of the decomposition of the exponential into its power series.
Here, in the spirit of Manicini--Sandeep--Tintarev's inequality on hyperbolic spaces,
we prove the following Moser--Trudinger inequality:
\begin{theorem}[A]
Let $(X,d,\mu)$ be an infinite volume $\rm{CD}(k,n)$-space, with nonzero Cheeger isoperimetric constant, such that
\begin{equation}
	\calV(X)=\inf_{x\in X}\mu(B_x(1))>0\,.
\end{equation}
	Denote by $\Phi$ the isoperimetric profile of $(X,d,\mu)$.
	Let $m\geq 2$ be an integer such that
	\begin{equation}
		l_\infty^{(m)}=\liminf_{t\rightarrow 0}\frac{1}{\omega_m}\frac{\Phi(t)^m}{m^mt^{m-1}}>0\,.
	\end{equation}
	Then for any $\alpha\leq m(l^{(m)}_\infty\omega_{m-1})^\frac{1}{m-1}$, there is a constant $C>0$ such that for any
	$u\in W^{1,m}(X)$ satisfying $\|\rm{lip}(u)\|_m\leq 1$, we have
\begin{equation}
	\int_X \big(\exp\big(\alpha u^\frac{m}{m-1}\big)-P_m(\alpha u^\frac{m}{m-1})\big)d\mu\leq C\,.
\end{equation}
\end{theorem}
Note that $l^{(n)}_\infty(X)$ is always nonzero for a $\rm{CD}(k,n)$-space, so applying this theorem with $m=n$ yields the classical Moser--Trudinger inequality for Riemannian manifolds.
See Section 4 for the proof of this theorem.
This theorem characterizes the manifolds admitting that Moser--Trudinger inequality among
manifolds with Ricci curvature bounded below. The presented proof relies on a rearrangement method,
allowing us to compare with the model case of the hyperbolic space with a conical singularity,
that can be handled as a consequence of \cite{MST13}. Note that this theorem extends the results of Yang obtained in \cite{Yan12} for asymptotically negatively curved manifolds.
Rearrangement methods and Polya--Szegö inequalities in metric measure spaces have been recently successfully used to derive other functional inequalities, see \cite{MS20,NV24}.

It is quite standard from such an inequality to derive log-Sobolev embeddings for metric measure spaces, which have been the focus of recent interest, see Balogh--Kristály--Tripaldi~\cite{BKT24}.

The proof presented by rearrangements is also  adapted to the study of compact manifolds,
or the study of compact $\rm{CD}(k,n)$-spaces. In this setup, we prove a Moser--Trudinger type inequality, specifically for zero-median functions:
\begin{theorem}[B]
	Let $(X,d,\mu)$ be a compact $\rm{CD}(k,n)$ space.
	Denote by $\Phi$ the isoperimetric profile of $(X,d,\mu)$.
	Let $m\geq 2$ be an integer such that
	\begin{equation}
		l_\infty^{(m)}=\liminf_{t\rightarrow 0}\frac{1}{\omega_m}\frac{\Phi(t)^m}{m^mt^{m-1}}>0\,.
	\end{equation}
	Then for any $\alpha\leq m(l^{(m)}_\infty\omega_{m-1})^\frac{1}{m-1}$, there is a constant $C>0$ such that for any
	$u\in W^{1,m}(X)$ satisfying $\|\rm{lip}(u)\|_m\leq 1$ with median $c$, we have
	\begin{equation}
		\int_X \big(e^{\alpha (u-c)^\frac{m}{m-1}}-P_m(\alpha (u-c)^\frac{m}{m-1})\big)d\mu\leq C\,.
	\end{equation}
	If $(X,d,\mu)$ is a smooth compact manifold $(M,g)$,
	the constant $C(\alpha)$ can be replaced by a constant depending only on
	the range of the curvature, the injectivity radius and the Cheeger constant of $M$.
\end{theorem}
And another one, in the same spirit, for zero-average functions:
\begin{theorem}[C]
	Let $(X,d,\mu)$ be a compact $\rm{CD}(k,n)$-space whose volume is larger than $V_0>0$
	Denote by $\Phi$ the isoperimetric profile of $(X,d,\mu)$.
	Let $m\geq 2$ be an integer such that
	\begin{equation}
		l_\infty^{(m)}=\liminf_{t\rightarrow 0}\frac{1}{\omega_m}\frac{\Phi(t)^m}{m^mt^{m-1}}>0\,.
	\end{equation}
	Then for any $\alpha\leq m(l^{(m)}_\infty\omega_{m-1})^\frac{1}{m-1}$, there is a constant $C>0$ such that for any
	$u\in W^{1,m}(X)$ satisfying $\|\rm{lip}(u)\|_m\leq 1$ with zero average, we have
	\begin{equation*}
		\int_X\big(e^{\alpha u^\frac{m}{m-1}}
		-P_m(\alpha u^\frac{m}{m-1})\big)d\mu \leq C(X)\,.
	\end{equation*}
	If $(X,d,\mu)$ is a smooth compact manifold $(M,g)$,
	the constant $C(X)$ can be replaced by a constant depending only on $V_0$,
	the range of the curvature, the injectivity radius, and the Cheeger constant.
\end{theorem}
See Section 5 for the proof of both these theorems.
As far as the author is concerned, the dependency of the constant involved here is new, especially
the fact that it can be made Volume independent. It seems in general to be a hard problem to
find the explicit value of $C$ appearing here, a question raised by Dolbeault--Esteban--Jankowiak in \cite{DEJ14}. The proof presented here relying on rearrangements, we don't expect the constant $C$ found to be sharp.
However it allows to control $C$ in large volumes, and the author used a version of that fact for hyperbolic surfaces in \cite{Bro23H4} to study minimal surfaces in the hyperbolic
$4$-space, and construct examples of hyperbolic structures on nontrivial disc bundles over
a surface.
\subsection{Acknowledgements}
The author thanks Nicolas Tholozan and Paul Laurain for their help and discussions on that topic a few years ago.
The author would also like to thank Daniele Semola and Gioacchino Antonelli
for their precious mathematical help on this topic.
The author thanks the Max Planck Institute for Mathematics in the Sciences (Leipzig) for the quality of its working environment.
This project has received funding from the European Research Council (ERC) under the European Union’s Horizon 2020 research and innovation programme (grant agreement No 101018839), ERC GA 101018839.
The author acknowledges support of the Institut Henri Poincaré (UAR 839 CNRS-Sorbonne Université) and LabEx Carmin (ANR-10-LABX-S9.09).

\section{Context, geometric preliminaries}
Here we introduce the geometric notions needed for rearrangements and for the study of the Moser--Trudinger inequality.
\subsection{CD, RCD, and ncRCD(k,n)-spaces and manifolds}
When generalizing analysis on manifolds to metric measure spaces, various "well-behaved" families appear in the literature.
We give here a quick introduction on these concepts.
\begin{definition}
	A \emph{metric measure space} is a a triple $(X,d,\mu)$, where
	$(X,d)$ is a complete, separable metric space and $\mu$ is a Borel measure on $X$, that is positive, finite on balls and not the zero measure.
\end{definition}
In a metric measure space $(X,d,\mu)$, for any $x\in X$ the function $r\mapsto \mu(B(x,r))$ is a monotonous function. As such, it is almost everywhere differentiable,
meaning that the perimeter of the sphere $S(x,r)$ is finite for almost every $r>0$. To simplify the statements to come, we always assume that the support of $\mu$ is $X$ itself.
We work with the classical notion of Perimeter in a metric measure space,
\begin{definition}
	Let $(X,d,\mu)$ be a metric measure space. Let $A\subset X$ measurable set.
	For $f:X\rightarrow\R$ and  $x\in X$, we define the slope of $f$ at $x$:
	\begin{equation}
		\rm{lip}(f)(x)=\limsup_{y\rightarrow x}\frac{|f(x)-f(y)|}{d(x,y)}\,.
	\end{equation}
	Denote by $\rm{Lip}(X)$ the space of Lipschitz functions on $X$, those whose slope is uniformly bounded.

	We define the Perimeter of $A$, denoted $\rm{Per}(A)$, the following quantity:
	\begin{equation}
		\rm{Per(A)}=\inf\bigg\{\liminf\int_X\rm{lip}(f_n) d\mu:\,f_i\in\rm{Lip}_{\rm{loc}}(X),\,f_i\rightarrow \chi_A\text{ in }L^1_{\rm{loc}}(X)\bigg\}\,.
	\end{equation}
\end{definition}
The completeness and separability condition is enough to satisfy a first essential theorem in Analysis: the Co-area formula.
In a metric measure space, the Sobolev spaces are classically defined in the following way (see Hajlasz~\cite{Haj96})
\begin{definition}
	Let $(X,d,\mu)$ be a metric measure space. Let $p\geq 1$.
	For $f:X\rightarrow\R$ measurable, its $p$-Cheeger energy, denoted $\rm{Ch}_p(f)$, is the following quantity:
	\begin{equation}
		\rm{Ch}_p(f)=\int_X \rm{lip}(f)^pd\mu\,.
	\end{equation}
	The Sobolev space $W^{1,p}(X)$ is the following Banach vector space:
	\begin{equation}
		W^{1,p}(X)=\big\{f\in L^p(X,\mu):\,\rm{Ch}_p(f)<\infty\big\}\,,
	\end{equation}
	endowed with the norm
	\begin{equation}
		\|f\|_{1,p}=(\|f\|^p_p+\rm{Ch}_p(f))^\frac{1}{p}\,.
	\end{equation}
\end{definition}
Sturm~\cite{Stu06I} introduced a condition called Curvature-Dimension (CD) condition, that is enough to satisfy Bishop--Gromov's inequality. As such, they are natural candidates for a generalization of manifolds with Ricci curvature bounded below.
The statement of that condition relies on the notion of relative entropy between measures, that is the following:
\begin{definition}
	Let $(X,d)$ be a metric space. We denote by $\mathcal{P}_2(X)$ the space of probability measures on $X$ with finite $L^2$-moment, i.e. for an (arbitrary) choice of basepoint $o$,
	\begin{equation}
		\int_X d(x,o)^2 d\mu(x)<\infty\,.
	\end{equation}
	The Wasserstein metric on $\calP_2(X)$ is the following:
	\begin{equation}
		d_W(\mu,\nu)=\inf\big\{\big(\int_{X\times X} d^2(x,y) dq(x,y)\big)^2:\,q\text{ coupling of }\mu,\nu\big\}
	\end{equation}
	where a coupling is a measure $q$ on $X\times X$ satisfying $q(A\times X)=\mu(A)$, $q(X\times A)=\nu(A)$ for any measurable set $A\subset X$.
	The relative entropy of two positive measures $\mu,\nu$ on $X$ is the following quantity:
	\begin{equation}
		\rm{Ent}(\nu|\mu)=\int_X \frac{d\nu}{d\mu}\log\big(\frac{d\nu}{d\mu}\big)d\mu\,.
	\end{equation}
\end{definition}
Sturm proved that $(\calP_2(X),d_W)$ is a complete separable metric space  containing $(X,d)$ totally geodesically.
\begin{definition}[$CD(K,\infty)$-spaces]
	Let $K\in\R$.
	A metric measured space $(X,d,\mu)$ satisfies the $CD(K,\infty)$ condition if for any $\nu_0,\nu_1\in\calP_2(X)$ with finite relative entropy to $\mu$,
	there exists a $d_W$-geodesic $\Gamma:[0,1]\rightarrow\calP_2(X)$ from $\nu_0$ to $\nu_1$ with:
	\begin{equation}
		\rm{Ent}(\Gamma(t)|\mu)\leq(1-t)\rm{Ent}(\nu_0|\mu)+(1-t)\rm{Ent}(\nu_1|\mu)-\frac{K}{2}t(1-t)d^2_W(\nu_0,\nu_1),\,\forall t\in [0,1]\,.
	\end{equation}
\end{definition}
The definition of $CD(K,N)$-space relies on a bound on another kind of Entropy, which is more subtle to write. We refer the reader to Sturm~\cite{Stu06II}, and Lott--Villani~\cite{LV07} for a definition. Note that Sturm's definition is a priori weaker, one gets to Lott--Villani's definition by adding an "essentially non-branching condition", which is essential
for most results presented here. For the remainder of the paper, by $\rm{C}(k,n)$-space we will mean an essential non-branching $\rm{CD}(k,n)$-space,
suitable in particular to use Rajala's results~\cite{Raj12}.
The main fact we will use about $CD(K,N)$-spaces will be the validity of Bishop--Gromov's comparison theorem, cf Sturm~\cite{Stu06II}, Theorem 2.3., Villani~\cite{Vil09}, Theorem 30.11.:
\begin{theorem}[Bishop--Gromov's Volume comparison theorem]
Let $X$ be a~$\rm{CD}(k,n)$-space, with $n>1$. Let $v(n,k,r)$ denote the volume of the ball of radius $r$. Let $T=\infty$ if $k\leq 0$ and $T=\pi\sqrt{(n-1)/k}$ if $k>0$.
in the simply connected Riemannian manifold of constant sectional curvature~$k$ and dimension~$n$,
and by $s(n,k,r)$ the area of the sphere of radius $r$ in it.
Then for any $x\in X$, the function $\frac{\Vol(B_r(x))}{v(n,k,r)}$ is non-increasing in $r$,
and for every $r\leq R\leq T$,
\begin{equation}
	\frac{\rm{Per}(B_R(x))}{s(n,k,R)}\leq\frac{\rm{Per}(B_r(x))}{s(n,k,r)}\,.
\end{equation}
Moreover, for any $0<r<T$,
\begin{equation}
	\frac{\rm{Per}(B_r(x))}{s(n,k,r)}\leq\frac{\Vol(B_r(x))}{v(n,k,r)}\,.
\end{equation}
\end{theorem}
Remark that this last inequality is a straightforward consequence of the monotonicity property from~\cite{Vil09}. Also, it implies that in $CD(k,n)$-spaces, the measure has no atoms, and is actually zero on all spheres.
In $CD(k,n)$-spaces, some functional inequalities are already true, like Poincaré inequalities.

Introduced by Gigli in~\cite{Gig15}, RCD(k,n)-spaces are a family of metric measure spaces
satisfying the $CD(k,n)$-curvature-dimension inequality, and infinitesimally Hilbertian, i.e. its Cheeger Energy is a quadratic form:
\begin{definition}
	Let $(X,d,\mu)$ be a metric measure space.
	We say that $(X,d,\mu)$ is infinitesimally Hilbertian if its $2$-Cheeger energy $\rm{Ch}$ is a quadratic form on $L^2(\mu)$, i.e. if $W^{1,2}(X)$ is a Hilbert space.
	We say that $(X,d,\mu)$ is a $\rm{RCD}(k,n)$-space if it is a $\rm{CD}(k,n)$-space and it is infinitesimally Hilbertian.
\end{definition}
We continue with some regularity assumptions that will be used later. The definition of noncollapsed $\rm{RCD}(k,n)$-space, or $\rm{ncRCD}(k,n)$-space, comes from De Philippis--Gigli~\cite{DG18}.
\begin{definition}
	An $\rm{RCD}(k,n)$-space $(X,d,\mu)$ is \emph{noncollapsed} if the measure $\mu$ coincides with the $n$-dimensional Hausdorff measure $\calH^n$.
\end{definition}
De Philippis and Gigli studied the behavior of the volume of small balls in those spaces. They considered the Bishop--Gromov density:
\begin{definition}
	Let $(X,d,\mu)$ be a metric measure space. Let $x\in X$
	The \emph{Bishop--Gromov density} at $x$ is the following quantity:
	\begin{equation}
		\theta(x)=\lim_{r\rightarrow 0}\frac{\mu(B(x,r))}{\omega_n r^n}\,.
	\end{equation}
	where $\omega_n=\frac{\pi^{n/2}}{\int_0^\infty t^{n/2} e^{-t}dt}$ is the volume of the unit ball in $\R^n$ (if $n$ is an integer).
\end{definition}
\begin{proposition}
	Let $(X,d,\calH^n)$ be a non-collapsed $\rm{RCD}(k,n)$-space. Then for $\calH^n$-almost every $x\in X$,
	\begin{equation}
		\theta(x)<\infty\,.
	\end{equation}
\end{proposition}
Finally, a last assumption is crucial when studying the isoperimetric problem on a metric measure space. It appeared already when doing analysis on manifolds, see Hebey~\cite{Heb96}.
\begin{definition}
	Let $(X,d,\mu)$ be a metric measure space. We denote by $\calV(X)$ the infimum of the volume of unit balls in $X$:
	\begin{equation}
		\calV(X)=\inf_{x\in X}\mu(B(x,1))\,.
	\end{equation}
\end{definition}
\begin{lemma}\label{lem:cpct}
	Let $(X,d,\mu)$ be a compact metric measure space. Then $\calV(X)>0$.
\end{lemma}
\begin{proof}
	This relies on the assumption that $\mu$ has full support. Indeed, assume that $\calV(X)=0$.
	Then there is a sequence $(x_n)$ in $X$ such that $\mu(B_{x_n}(1))\rightarrow 0$. Since $X$ is compact, up to extraction $x_n\rightarrow x\in X$.
	But then for large enough $n$, $B_x(\frac{1}{2})\subset B_{x_n}(1)$.
	Hence we deduce that $\mu(B_x(\frac{1}{2})=0$, which is absurd since $\mu$ has full support. So we have shown that $\calV(X)>0$, as claimed.
\end{proof}
\subsection{The Coarea formula for metric measure spaces}
Here we set the Coarea formula we need in order to make rearrangements on metric measure spaces. See Mondino--Semola~\cite{MS20} for a more detailed description.
The natural functions on which to apply the Coarea formula are the bounded variation functions:
\begin{definition}
	Let $(X,d,\mu)$ be a metric measure space. A measurable function $u:X\rightarrow\R$ is of bounded variation if there exists a sequence of locally Lipschitz functions $f_n$ converging to $f$ in $L^1(\mu)$
	satisfying
	\begin{equation}
		\limsup_{n\rightarrow+\infty}\int \rm{lip(f_n)}d\mu<+\infty\,.
	\end{equation}
	We denote by $\rm{BV}(X)$ or $\rm{BV}(X,\mu)$ the space of bounded variation functions on $X$.
	For any $f\in\rm{BV}(X)$, we denote by $\|Df\|_\ast$ the Borel measure such that for every open set $A$,
	\begin{equation}
		\|Df\|_\ast(A)=\inf\big\{\liminf_{n\rightarrow\infty}\int_A\rm{lip}(f_n)d\mu,\,f_n\in\rm{Lip}_{\rm{loc}}(A), f_n\rightarrow f\text{ in } L^1(A)\big\}\,.
	\end{equation}
	For $E\subset X$ a set such that $\mathbf 1_E$ is of bounded variation, we denote by $\|\partial E\|=\|D\mathbf 1_E\|$.
\end{definition}
Remark that by construction, $\|\partial E\|(X)=\rm{Per}(E)$. Also, for any $f\in\rm{Lip}_{\rm{loc}}(X)$, $\|Df\|\leq\rm{lip}(f)\mu$.
The following coarea formula for metric measure spaces is due to Miranda~\cite{Mir03}:
\begin{theorem}
	Let $(X,d,\mu)$ be a metric measure space, $u\in\rm{BV}(X)$, $v:X\rightarrow\R$ measurable and $A$ a Borel set in $X$.
	Denote by $E_t=\{u>t\}$ the surlevel sets.
	Then
	\begin{equation}
		\int_\R\big(\int_A v(x) d\|\partial E_t\|(x)\big)dt=\int_A v(x)d\|Du\|(x)\,.
	\end{equation}
\end{theorem}
Note that it is a consequence of the theorem that $E_t$ is of finite perimeter for almost every $t$.
We will also make use of the following version of the coarea formula:
\begin{theorem}\label{thm:coareastrong}
	Let $(X,d,\mu)$ be a metric measure space. Let $\Omega\subset X$ be a finite volume open subset.
	Let $u\in\rm{Lip}_{\rm{loc}}(\Omega)$ be non-negative, and satisfying $\frac{d\|D u\|_1}{d\mu}(x)\neq 0$ for $\mu$-almost every $x\in\{u>0\}$.
	Then its distribution $A:\R^\ast_+\rightarrow[0,\mu(\Omega)]$ is absolutely continuous and verifies, for any $t_1<t_2$:
	\begin{equation}
		A(t_1)-A(t_2)=\int_{t_1}^{t_2}\big(\int(\frac{d\|Du\|_1}{d\mu})^{-1}d\|\partial E_t\|\big)dt\,.
	\end{equation}
	At the limit when $t_2\rightarrow t_1=t$ we obtain for almost every $t$,
	\begin{equation}
		-A'(t)=\int(\frac{d\|Du\|}{d\mu})^{-1}d\|\partial E_t\|\,.
	\end{equation}	
\end{theorem}
Finally, we want to approximate any function in a Sobolev space by a locally Lipschitz function satisfying this non-vanishing of the gradient condition.
Mondino--Semola proved that this is possible provided the metric measure space is locally compact and geodesic. This is their Lemma 3.6.
\begin{proposition}\label{prop:approx}
	Let $(X,d,\mu)$ be a metric measure space which is locally compact, geodesic. Let $\Omega\subset X$ be a finite volume open set.
	For any non-negative $u\in\rm{Lip}_c(X)$ with $\int\rm{lip}(u)^pd\mu<\infty$, there exists a sequence $(u_n)$ compactly supported, Lipschitz, non-negative,
	such that $\frac{d\|Du_n\|}{d\mu}\neq 0$ for $\mu$-almost every $x\in\{u_n>0\}$, and such that
	\begin{equation}
		u_n\rightarrow u\text{ in }W^{1,p}(X,\mu)\,.
	\end{equation}
\end{proposition}
Note that locally compact geodesic metric measure spaces contains $\rm{CD}(k.n)$-spaces, as well as weighted Riemannian manifolds $(M,g,\mu)$.
\subsection{The Cheeger constant and Buser-type inequalities}
Let $(X,d,\mu)$ be a metric measure space. Its Cheeger constant is defined in the following way:
\begin{definition}[Cheeger constant]\label{def:cheeger}
The \emph{Cheeger constant} of $X$, denoted $h(X)$,
is the following constant:
\begin{equation}
	h(X)=\inf\{\frac{\rm{Per}(A)}{\mu(A)},\,A\subset X, \overline A\text{ compact},\,
	2\Vol(A)\leq\Vol(X)\}
\end{equation}
\end{definition}
The Cheeger quotient is heavily related to the $p$-th spectral gap,
which we define here for $p>1$.
\begin{definition}
	Let $(X,d,\mu)$ be a metric measure space.
	If $\Vol(X)=\infty$, we call $p$-th spectral gap the quantity
	\begin{equation}
		\lambda_p=\inf\big\{\frac{\rm{Ch}_p(f)}{\int_X |f|^p},\,
			f\in W^{1,p}(X)\big\}\,.
	\end{equation}
	If $\Vol(X)<\infty$, we call $p$-th spectral gap the quantity
	\begin{equation}
		\lambda_p=\inf\big\{\frac{\rm{Ch}_p(f)}
		{\inf_{c\in\R}\int_X |f-c|^p},\, f\in W^{1,p}(X)\big\}\,.
	\end{equation}
\end{definition}
Note that $\lambda_2$ is the usual spectral gap of a manifold. In general,
$\lambda_p$ is a synthetic spectral gap for the $p$-Laplacian operator on $W^{1,p}(M)$, when it exists.
The celebrated Cheeger inequality gives a lower bound on $\lambda_2$ in terms of the Cheeger constant. It holds in metric measure spaces (see De Ponti--Mondino~\cite{DM21}, Theorem 4.2.):
\begin{theorem}
	Let $(X,d,\mu)$ be a metric measure space. Then its spectral gap $\lambda_2$ and its Cheeger constant $h$ satisfy:
	\begin{equation}
		\lambda_p\geq\frac{1}{p^p}h^p\,.
	\end{equation}
\end{theorem}
\begin{proof}
	The proof in the case $p=2$ is exactly the one from De Ponti--Mondino~\cite{DM21}. We write it here for arbitrary $p>1$.
	Let $\eps>0$. By construction, there is $f\in\rm{Lip}(X)$ such that
	\begin{equation*}
		\lambda_p\geq \frac{\int_X\rm{lip}(f)^pd\mu}{\int_X|f|^pd\mu}-\eps\,.
	\end{equation*}
	Bounding the slope of $|f|^p$ in the following way:
	\begin{equation*}
		\rm{lip}(|f|^p)\leq p|f|^{p-1}\rm{lip}(f)\,,
	\end{equation*}
	we obtain, vie a Hölder inequality:
	\begin{equation*}
		\int_X\rm{lip}(|f|^p)d\mu\leq p\big(\int_X|f|^pd\mu\big)^\frac{p-1}{p}\big(\int_X\rm{lip}(f)^pd\mu\big)^\frac{1}{p}\,.
	\end{equation*}
	Combining this with the above inequality, we obtain
	\begin{equation*}
		\lambda_p\geq\frac{\big(\int_X\rm{lip}(|f|^p)d\mu\big)^p}{p^p\big(\int_X|f|^pd\mu\big)^p}-\eps\,.
	\end{equation*}
	Now to control the right-hand side we use the co-area formula:
	\begin{equation*}
		\int_X\rm{lip}(|f|^p)d\mu=\int_0^{\sup|f|^p}\rm{Per}(\{|f|^p>t\})dt\,.
	\end{equation*}
	From which we deduce that
	\begin{flalign*}
		\int_X\rm{lip}(|f|^p)d\mu&\geq h\int_0^{\sup|f|^p}\mu(\{|f|^p>t\})dt\\
			&\geq \int_X |f|^pd\mu\,.
	\end{flalign*}
	All in all, we have shown that
	\begin{equation*}
		\lambda_p\geq\frac{1}{p^p}h^p-\eps\,.
	\end{equation*}
	for any $\eps>0$, as desired.
\end{proof}
For manifolds with Ricci curvature bounded below, Matei~\cite{Mat00}, generalizing on the arguments of Buser~\cite{Bus82} proved an upper bound on the $p$-th
spectral gap from the Cheeger constant.
\begin{theorem}\label{thm:buser}
	Assume $(M,g)$ is a noncollapsed
	manifold with Ricci curvature bounded below.
	Then there is a constant $C(k,n,p)$ such that
	\begin{equation}
		\lambda_p\leq c(n,k,p)(h+h^p)\,.
	\end{equation}
\end{theorem}
To the best of our knowledge, this inequality has not been proved in
the general setup of $\rm{CD}(k,n)$-spaces yet, apart from the case $p=2$ for $RCD(k,n)$-spaces, see De Ponti--Mondino~\cite{DM21}.

A control on the Cheeger constant will allow us to control the Moser--Trudinger inequality on the large scale. At small scale, a quantity that will play
a critical role will be the asymptotic growth of balls.
\subsection{The asymptotic growth ratio}
\begin{definition}\label{def:agr}
	Let $(X,d,\mu)$ be a metric measure space. Let $n\geq 1$
	Then its $n$-asymptotic growth ratio is the quantity $\theta_\infty$:
	\begin{equation}
		\theta_\infty=\liminf_{r\rightarrow 0}\inf_{x\in X}
		\frac{\Vol(B(x,r))}{\omega_n r^n}\,.
	\end{equation}
\end{definition}
Note that on an $n$-dimensional manifold, the $n$-asymptotic growth ratio is always less than $1$.
Unless otherwise specified, the asymptotic growth ratio of a $\rm{CD}(k,n)$-space always means its $n$-asymptotic growth ratio.
Standard considerations on the volume growth of balls shows that it is $1$ for any smooth compact manifold, or on a manifold with bounded sectional curvature and lower bounded injectivity radius.
Antonelli et.al.\cite{APP+24} have studied this ratio, and shown that it is nonzero for noncollapsed $\rm{RCD}(k,n)$-spaces.

\section{Domination by radial spaces}
\subsection{Radial spaces and the isoperimetric profile}
Let $n\geq 2$ and $\E^n$ denote the real vector space of dimension $n$.
\begin{definition}
A metric $\sigma$ on $\E^n$ is said to be \emph{radial} if there is a continuous
function~$g:\R_+\rightarrow\R_+^\ast$ such that $\sigma$ can be written, in polar coordinates:
\begin{equation}
\sigma=dr^2+g(r)^2d\theta^2\,,
\end{equation}
	where $d\theta^2$ denotes the round metric on the unit sphere $\S^{n-1}$.
\end{definition}
Note that if we want to avoid a cone-type singularity at zero,
the function $g$ must satisfy $\frac{g(r)}{r}\rightarrow 1$ when $r$ goes to zero.
Since the function $g$ is chosen continuous, $(\E^n,\sigma)$ is a Riemannian manifold, and such the notion of length measure and gradient are well-defined and do not
need the the setup of metric measure spaces.

A convenient notion is that of the radial profile of such a space:
\begin{definition}
Let $(\E^n,\sigma)$ be a radial space. Then its \emph{radial profile} $\varphi_\sigma$ is the function which
to the centered ball of volume $t$ gives its perimeter. More precisely,
	denoting by~$\calH^{n-1}$
	the $(n-1)$-dimensional Hausdorff measure on $(\E,\sigma)$,
	the radial profile is the function
\begin{equation}
\varphi_\sigma:\left\{\begin{array}{ll}
	[0,\Vol(\E^n,\sigma)]&\mapsto \R_+^\ast\\
	t&\rightarrow \calH^{n-1}(S(0,T))\end{array}\right.
\end{equation}
where $T$ is the unique $T>0$ such that $\Vol(B(0,T))=t$.
\end{definition}
Note that $T$ is well-defined since we assumed $g>0$. Also, the function $T\mapsto \Vol(B(0,T))$
is $C^1$-regular, and an implicit function theorem shows that $\varphi_\sigma$ is continuous,
strictly increasing if $g$ is.
\begin{definition}
	Let $0<\alpha<1$. We call \emph{hyperbolic trumpet of cone angle $\alpha$},
	denoted $\H^n_\alpha$, the radial space with the radial metric
	$\E^n,dr^2+\alpha^\frac{2}{n-1}\sinh^{2}(r)d\theta^2$.
\end{definition}
\begin{remark}
	The hyperbolic trumpet of cone angle $\alpha$ is the pullback of the hyperbolic metric on $\E^n$ by the map $(r,\theta)\mapsto (r^\alpha,\theta)$, hence the cone angle singularity at zero.
	Note that it is not locally hyperbolic, even if its large scale geometry is similar.
\end{remark}
\begin{proposition}
	Let $f:\R_+\rightarrow\R_+$ be a continuous function vanishing only at zero,
	and such that $\frac{f(t)}{t^{1-\frac{1}{n}}}\rightarrow \alpha>0$ at zero.
	Then there is a radial space with metric $dr^2+g(r)^2d\theta$ whose radial profile
	is $f$. It admits a cone singularity of angle $\frac{\alpha^n}{n^{n-1}\omega_{n-1}}$ at zero.

	If moreover, there is $h>0$ and $t_0>0$ such that $f(t)=ht$ for $t\geq t_0$,
	then the metric is radially bilipschitz conformal to a hyperbolic plane with cone singularity $\frac{\alpha^n}{n^{n-1}\omega_{n-1}}$.
\end{proposition}
\begin{proof}
	For $g$ a $C^0$-regular function, $f$ is the radial profile of $dr^2+g(r)^2d\theta^2$
	if and only if:
	\begin{equation*}
		\forall s>0,\, G'(s)=\frac{1}{\omega_{n-1}}f(\omega_{n-1} G(s))\,,\, G(s)=\int_0^s g(u)^{n-1}du\,.
	\end{equation*}
	In order to solve it locally, introduce the function $F(s)=G(s)^\frac{1}{n}$.
	Then it must satisfy the ODE:
	\begin{equation*}
		F'(s)=\frac{1}{n\omega_{n-1}}\frac{f(\omega_{n-1} G(s))}{G(s)^{1-\frac{1}{n}}}
		=\frac{1}{n\omega_{n-1}}\frac{f(\omega_{n-1}F(s)^n)}{F(s)^{n-1}}\,.
	\end{equation*}
	From the asymptotics of $f$ at zero, we can see that there is a unique local
	solution $F$ to this ODE such that
	\begin{equation*}
		F(0)=0,\, F'(0)=\frac{\alpha}{n\omega_{n-1}^{\frac{1}{n}}}>0\,.
	\end{equation*}
	Since $f$ is positive continuous, it is clear that this solution can be extended uniquely to a solution defined on $\R_+$, nondecreasing, vanishing only at zero.

	This leads to corresponding functions $G$ and $g$, satisfying the asymptotics
	\begin{equation*}
		\frac{g(s)}{s}\underset{s\rightarrow 0}{\rightarrow}\frac{\alpha^\frac{n}{n-1}}{n\omega_{n-1}^\frac{1}{n-1}}\,.
	\end{equation*}
	Hence the metric is nonsingular at zero if and only if $\alpha=n^\frac{n-1}{n}\omega_{n-1}^\frac{1}{n}$, and otherwise it admits a cone type singularity, of total angle $\frac{\alpha^n}{n^{n-1}\omega_{n-1}}$.

	Assume now that there exists $h>0$ such that $f(t)= ht$ for $t\geq t_0$.
	That means that for $t\geq t_0$ we have
	\begin{equation*}
		F(t)=F(t_0)e^{\frac{h}{n}(t-t_0)},\, G(t)=G(t_0)e^{h(t-t_0)},\, g(t)=g(t_0)e^{\frac{h}{n-1}(t-t_0)}\,.
	\end{equation*}
	In particular the identity map is a Bilipschitz conformal transform
	between $\B^n, dr^2+g(r)^2d\theta^2$ and $\B^n, \big(\frac{\alpha^\frac{n}{n-1}}{hn\omega_{n-1}^\frac{1}{n-1}}\big)^2\sinh(hr)^2d\theta^2$. Composing with the radial map $h\mapsto hr$, we get a bilipschitz
	conformal transform to a hyperbolic space with cone angle $\frac{\alpha^n}{n^{n-1}\omega_{n-1}}$, as desired.
\end{proof}
Given a metric measure space $(X,d,\mu)$, we consider its isoperimetric profile:
\begin{definition}
	The isoperimetric profile $\Phi$ of $(X,d,\mu)$ is defined for $0\leq t\leq \mu(X)$:
\begin{equation}
	\Phi(t)=\inf\{\rm{Per}(A),\,A\subset X,\,\mu(A)=t, \overline A \text{ compact }\}
\end{equation}
\end{definition}
The isoperimetric profile has been the object of a lot of research, specifically in the field
of isoperimetric inequalities, see for instance \cite{Gal88,BM82}.
Here we will say that our manifold is dominated if its isoperimetric profile is larger
than the radial profile of a space of same dimension.
\begin{definition}
We say that $(X,d,\mu)$ is \emph{dominated} by the radial space $(\E^n,\sigma)$ if its isoperimetric profile
	satisfies $\Vol(\E^n)\geq\frac{\mu(X)}{2}$ and
\begin{equation}
	\forall t\leq\frac{\mu(X)}{2},\, \Phi(t)\geq\varphi_\sigma(t)\,.
\end{equation}
\end{definition}
Note that the definition still makes sense when $\mu(X)=\infty$ as well, in which case the
requirement $\sigma(\E^n)=\infty$ is equivalent to the completeness of $\sigma$.
We finish this subsection with the notion of "well-behaved" radial space, in which isoperimetric and
radial profile coïncide:
\begin{definition}
	We say that the radial space $(\E^n, dr^2+g(r)^2d\theta^2)$ is \emph{well-behaved}
	if its isoperimetric profile equals its radial profile, i.e. it is well-behaved if it
	is self-dominated.
\end{definition}
\begin{lemma}\label{lem:wellbehaved}
	For any $\alpha\leq 1$ the hyperbolic trumpet with cone angle $\alpha$ is well-behaved.
\end{lemma}
\begin{proof}
	It is well-known that balls are solution to the isoperimetric problem in hyperbolic space,
	see Schmidt~\cite{Sch43}. For a hyperbolic trumpet of cone angle $\alpha$, pick $D$ a domain,
	denote by $\Vol(D)$ its volume and $\Vol_{hyp}(D)$ its volume for the hyperbolic metric.
	Similarly, consider $\calH^{n-1}(D)$ and $\calH^{n-1}_{hyp}(D)$ its perimeters.
	Than a quick computation shows that
	\begin{equation*}
		\Vol(D)=\alpha\Vol_{hyp}(D),\,\calH^{n-1}(D)\geq\alpha\calH_{hyp}^{n-1}(D)
	\end{equation*}
	with equality when $D$ is a centered ball.
	This proves that the centered balls are solutions to the isoperimetric problem
	in hyperbolic trumpets, hence their isoperimetric profile equals their radial profile.
\end{proof}
\subsection{The isoperimetric profile at small volumes}
Here we present a control on the isoperimetric profile at small volumes, generalizing Berard--Meyer~\cite{BM82} to metric measure spaces.
This generalization was done by Antonelli--Pasqualetto--Pozzetta~\cite{APP22}, Prop~3.20.
\begin{theorem}\label{thm:smalliso}
	Let $(X,d,\mu)$ be a $\rm{CD}(k,n)$-space with $\calV(X)>0$.
	Denote by $\Phi$ its isoperimetric profile.
	Then there are $C,\eta>0$ such that
	\begin{equation}
		\mu(E)\leq\eta\,\Rightarrow\,\rm{Per}(E)\geq C\mu(E)^\frac{n-1}{n}\,.
	\end{equation}
	equivalently,
	\begin{equation}
		0\leq t\leq\eta\,\Rightarrow\,\Phi(t)\geq C t^\frac{n-1}{n}\,.
	\end{equation}
	The constants $C,\eta$ can be chosen to be dependent only in $k,n,\calV(X)$.
\end{theorem}
\begin{definition}\label{def:asymisolength}
	Let $(X,d,\mu)$ be a metric measure space, with $\Phi$ its isoperimetric profile.
	We call \emph{isoperimetric dimension} the quantity:
	\begin{equation}
		n_X=\inf\{n\geq 1:\,\liminf_{t\rightarrow 0}\frac{\Phi(t)}{t^{1-\frac{1}{n}}}>0\}\,.
	\end{equation}
	We call the \emph{asymptotic isoperimetric ratio} the quantity:
	\begin{equation}
		l_\infty(X)=\liminf_{t\rightarrow 0}\frac{\Phi(t)^{n_X}}{n_X^{n_X}\omega_{n_X} t^{n_X-1}}\,.
	\end{equation}
	For arbitrary $n$, we denote by $l_\infty^{(n)}(X)$ the quantity:
	\begin{equation}
		l_\infty^{(n)}(X)=\liminf_{t\rightarrow 0}\frac{\Phi(t)^n}{n^n\omega_n t^{n-1}}\,.
	\end{equation}
\end{definition}
\begin{remark}
	From the previous theorem, one sees that the isoperimetric dimension of a $\rm{CD}(k,n)$-space with $\calV(X)>0$ is always less or equal to $n$.
	It is clear that $l_\infty^{(n)}=+\infty$ for $n>n_X$ and $l_\infty^{(n)}=0$ for $n<n_X$. The author is not aware of criteria ensuring that $l_\infty(X)\in (0,+\infty)$.
\end{remark}
\begin{example}
	\begin{itemize}
		\item For a smooth compact riemannian $n$-manifold $(M^n,g)$, its isoperimetric dimension is $n$ and its asymptotic
			isoperimetric ratio is $1$.
		\item For a compact manifold $(M^n,g)$ with finitely many cone-type singularities of angles $(\theta_1,\ldots,\theta_n)$,
			its isoperimetric dimension is $n$ and its asymptotic isoperimetric ratio is $\min\{\theta_1,\ldots,\theta_n\}$.
		\item For a smooth compact weighted Riemannian manifold $(M^n,g,f\cdot\rm{Leb})$, the isoperimetric dimension is $n$ and the asymptotic isoperimetric ratio
			belongs to $(0,+\infty)$.
		\item the isoperimetric dimension of a hyperbolic cusp is $+\infty$.
	\end{itemize}
\end{example}
In \cite{APP+24}, they give a simpler understanding of these quantities for noncollapsed $\rm{RCD}(k,n)$-space.
\begin{proposition}
	Let $(X,d,\mu)$ be a noncollapsed $\rm{RCD}(k,n)$-space, with $\calV(X)>0$.
	Then its isoperimetric dimension is $n_X=n$, and its asymptotic isoperimetric ratio equals its asymptotic volume ratio:
	\begin{equation}
		\theta_\infty=l_\infty\,.
	\end{equation}
\end{proposition}

\subsection{Rearrangement to a dominating radial space}
Here we fix $(X,d,\mu)$ a metric measure space dominated by a radial space $(\E^n,\sigma)$.
We want to build a rearrangement, that is a map 
$W^{1,p}(X)_{>0}\rightarrow W^{1,p}(\E^n)$, mapping $f$ to $f^\ast$, satisfying
\begin{equation}
	\|\nabla f^\ast\|_p\leq\|\rm{lip}(f)\|_p,\,\sigma(\{f^\ast\geq t\})=\mu(\{f\geq t\})\,\forall t\in\R.
\end{equation}
In order for this to be possible, a first condition is that~$\mu(X)=+\infty$.
We will need a different appropriate rearrangement in the finite volume case,
considering only positive functions whose support ihas volume less than half of $\mu(X)$.
\begin{theorem}
Let $(X,d,\mu)$ be an infinite volume locally compact, geodesic metric measure space
dominated by a radial space $(\E^n,\sigma)$.
Consider the radial rearrangement, from $W^{1,p}(X)_{\geq 0}$ to $W^{1,p}(\E^n,\sigma)$
which to a nonnegative $u\in W^{1,p}(X)$ associates $\hat u\in W^{1,p}(\E^n)$
the unique radial decreasing left-continuous function satisfying
\begin{equation}
	\forall t>0,\,\mu(\{u>t\})=\sigma(\{\hat u> t\})\,.
\end{equation}
Then it satisfies
\begin{equation}
	\|\nabla \hat u\|_p\leq \|\rm{lip}(u)\|_p\,.
\end{equation}
\end{theorem}
\begin{proof}
	Let $u\in W^{1,p}(X)\cap\rm{BV}(X)$ be a positive function with bounded variation.
	First assume that it is compactly supported and with non-vanishing gradient, i.e. $\frac{d\|Du\|}{d\mu}(x)\neq 0$ for $\mu$-almost every $x\in\{u>0\}$.
	Denote by $A(t)$ the distribution function of $u$, i.e. $A(t)=\mu(\{u>t\})=\mu(E_t)$ for all $t\geq 0$
	Applying the co-area formula from Theorem~\ref{thm:coareastrong}, we have for every $t$ in the range of $u$:
	\begin{equation*}
		A(t)=\int_{s\geq t}\big(\int(\frac{d\|Du\|}{d\mu})^{-1}d\|\partial E_s\|\big)ds\,.
	\end{equation*},
	In particular, $A$ is differentiable almost everywhere, and satisfies:
	\begin{equation*}
		-A'(t)=\int(\frac{d\|Du\|}{d\mu})^{-1}d\|\partial E_t\|\,.
	\end{equation*} 
	Also, denoting $l(t)=\rm{Per}(E_t)$, $l$ is finite for almost every $t$, and via Hölder's inequality:
	\begin{equation*}
		l(t)=\|\partial E_t\|(X)\leq\big(\int(\frac{d\|Du\|}{d\mu})^{p-1}d\|\partial E_t\|\big)^\frac{1}{p}
			\big(\int(\frac{d\|Du\|}{d\mu})^{-1}d\|\partial E_t\|\big)^\frac{p-1}{p}\,.
	\end{equation*}
	From which we deduce that $(-A')^{1-p}l^p$ is integrable, with the bound:
	\begin{equation*}
		\int_{\R^+} (-A'(s))^{1-p}l(s)^p ds\leq \int(\frac{d\|Du\|}{d\mu})^pd\mu\leq\int\rm{lip}(u)^pd\mu\,,
	\end{equation*}
	where for the last inequality we used that the density of $\|Du\|$ with regard to $\mu$ is smaller than $\rm{lip}(u)$.
	Now consider the radial rearrangement $\hat u$ on $(\E^n,\sigma)$, with its length measure $l^\sigma$, and the corresponding functions $\hat A$, $\hat l$ defined by
	\begin{equation*}
		\hat A(t)=\sigma(\{\hat u>t\})=A(t)\,,\,\hat l(t)=l^\sigma(\{\hat u =t\})\,.
	\end{equation*}
	By construction of the radial profile, we have
	\begin{equation*}
		\hat l(t)=\varphi(A(t))\,.
	\end{equation*}
	And since $\hat u$ is radial, the Hölder inequality used above becomes an equality, and so is the identification between the gradients,
	and we obtain
	\begin{equation*}
		\|\nabla\hat u\|_p^p=\int_{\R_+}(-A'(s))^{1-p}\hat l(s)^pds=\int_{\R_+}(-A'(s))^{1-p}\varphi(A(s))^pds\,.
	\end{equation*}
	Finally, because of the domination assumption, we know that
	\begin{equation*}
		\varphi(A(s))\leq\Phi(A(s))\leq l(s)\,.
	\end{equation*}
	Hence we have shown that
	\begin{equation*}
		\|\nabla \hat u\|_p\leq\int_{\R_+}(-A'(s))^{1-p}l(s)^p ds\leq\int \rm{lip}(u)^pd\mu\,,
	\end{equation*}
	as desired.

	Let now $u\in W^{1,p}(X)$ a positive function, not necessarily of bounded variation.
	Consider $u_k\rightarrow u$ a regularizing sequence of bounded variation functions with nonvanishing gradients, converging towards $u$ in $W^{1,p}(M)$.
	This exists since $(X,d,\mu)$ is locally compact geodesic (see Proposition~\ref{prop:approx}).
	
	From the discussion above, the sequence $\hat u_k$ is bounded in $W^{1,p}_{rad}(\E^n)$,
	hence there is $v$ such that $\hat u_k\rightarrow v$ weakly, and since they're radial we
	can assume they converge uniformly on compact subsets of $\E^n-\{0\}$.
	We claim that $v=\hat u$. The proof is standard, and can be found in
	Brothers--Ziemer~\cite{BZ88}. We reproduce it here for the reader's convenience.
	Let $t>0$ be a value such that $\mu(\{u=t\})=0$. This is true for all but at most countably 
	$t>0$, since $u^p$ is integrable. Define the sets
	\begin{equation*}
		E_t^k=\{u_k>t\},\, E_t=\{u>t\},\, \hat E_t^k=\{\hat u_k>t\}, E_t^v=\{v>t\}\,.
	\end{equation*}
	From the convergence of $u_k$, it can be seen that the symmetric difference satisfies:
	\begin{equation*}
		\mu(E_t^k\Delta E_t^u)\rightarrow 0\,.
	\end{equation*}
	Indeed, for any $\eta>0$ we can control the volume
	\begin{equation*}
		\mu(E_t^k\Delta E_t^u)\leq \mu(\{|u-u_k|\geq\eta\})+\mu(t-\eta\leq u\leq t+\eta)\,
	\end{equation*}
	Fix $\eps>0$. Since $\mu(\{u=t\}=0)$ and $u$ is locally integrable, there is $\eta>0$ such that
	\begin{equation*}
		\mu(|u-t|\leq \eta)\leq\frac{\eps}{2}\,.
	\end{equation*}
	Also, by Markov's inequality
	\begin{equation*}
		\mu(\{|u-u_k|\geq\eta\})\leq\frac{\|u-u_k\|_p^p}{\eta^p}
	\end{equation*}
	which is smaller than $\frac{\eps}{2}$ for $k$ large enough.
	Hence we have proven the convergence, for all but countably many values of $t$:
	\begin{equation*}
		\mu(\{u_k>t\})\rightarrow\mu(\{u>t\})
	\end{equation*}
	Similarly, the uniform convergence on compact sets implies that for all but countably
	many values of $t$,
	\begin{equation*}
		\sigma(\{\hat u_k>t\}\rightarrow\sigma(\{v>t\})\,.
	\end{equation*}
	All in all, for almost every $t$,
	\begin{equation*}
		\sigma(v>t)=\mu(u>t)
	\end{equation*}
	and $v$ is radially symmetric decreasing. Hence $v=\hat u$, at all but countably many values.
	Since at those values we can still characterize the volume
	\begin{equation*}
		\sigma(v>t)=\lim_{\eta\rightarrow 0}\sigma(v>t+\eta)
	\end{equation*}
	we obtain that $v=\hat u$, as claimed.
	Finally, since $v\in W^{1,p}$ and we have the weak convergence
	in~$W^{1,p}$ of $\hat u_k$ towards $v$, we deduce that
	\begin{equation*}
		\|\nabla\hat u\|_p\leq\liminf\|\nabla\hat u_k\|_p\leq\liminf\|\rm{lip}(u_k)\|_p
			=\|\rm{lip}(u)\|_p\,.
	\end{equation*}
	So $\|\nabla\hat u\|_p\leq\|\rm{lip}(u)\|_p$, as desired.
\end{proof}
\begin{remark}
	This rather technical process is necessary since the rearrangement is not continuous.
	See Almgren--Lieb~\cite{AL89} for a discussion of the noncontinuity of rearrangement.
	It seems, from the works of Brothers and Ziemer, that one could extend the proof
	to include the case of $W^{1,1}(M)$, however our main interest is in $W^{1,n}(M)$,
	so we won't do it here.
\end{remark}
Since the proof only used the fact that $u_k$ converges towards $u$ weakly, we have the following
corollary
\begin{corollary}
	Let $(X,d,\mu)$ be a locally compact geodesic metric measure space dominated by a radial space $(\E^n,\sigma)$,
	and $p>1$.
	Then the radial decreasing rearrangement from $W^{1,p}(X)$ to $W^{1,p}(\E^n)$
	which to $u$ associates $\widehat{|u|}$ is weakly continuous.
\end{corollary}
\subsection{Characterization of dominated metric measure spaces}
Being dominated by a nonsingular radial space implies strong properties of the isoperimetric profile.
It is natural to expect a $\rm{CD}(k,n)$-space to be dominated by a space bilipschitzly conformal to hyperbolic space.
This is however not true in all generality, as will be seen from the characterization of dominated spaces among noncollapsed $\rm{RCD}(k,n)$-spaces.
\begin{lemma}\label{lem:ballvolume}
	Let $(X,d,\mu)$ be a metric measure space, dominated by a radial space $(\E^n,\sigma)$, that is either nonsingular at zero, or admit a finite angle cone-singularity at zero.
	Then for any $r>0$ it satisfies
	\begin{equation}
		\inf_{x\in X}\mu(B_x(r))>0\,.
	\end{equation}
	Moreover,
	\begin{equation}
		\liminf_{r\rightarrow 0}\inf_{x\in X}\frac{\mu(B_x(r))}{\omega_n r^n}>0\,,\quad
	\end{equation}
\end{lemma}
\begin{proof}
	First, since the radial space $(\E^n,\sigma)$ admits a nonzero cone-singularity at zero,  we get the asymptotic for small balls:
	\begin{equation*}
		\sigma(B(0,r))\approx\theta\omega_n r^n,\,l^\sigma(S(0,r))\approx n\theta\omega_n r^{n-1}\,.
	\end{equation*}
	for some $\theta>0$ the cone angle at zero.
	Eventually we obtain that the radial profile $\varphi^\sigma$ satisfies the asymptotic for small $t$:
	\begin{equation*}
		\varphi^\sigma(t)\approx n(\theta\omega_n)^\frac{1}{n} t^{1-\frac{1}{n}}\,.
	\end{equation*}
	Via the domination property, and denoting $\Phi$ the isoperimetric profile of $(X,d,\mu)$, $\Phi\geq\varphi^\sigma$, hence there is $C>0$ and $\eta>0$ such that:
	\begin{equation*}
		0\leq t<\eta\quad\Rightarrow\quad \Phi(t)\geq C t^{1-\frac{1}{n}}\,.
	\end{equation*}
	Now fix $x\in X$, and $r_0>0$. If $\mu(B_x(r_0))<\eta$, for any $r\leq r_0$ we have
	\begin{equation*}
		\partial_r\mu(B_x(r))\geq\rm{Per}(B_x(r))\geq C\mu(B_x(r))^{1-\frac{1}{n}}\,.
	\end{equation*}
	From which we deduce the lower bound:
	\begin{equation*}
		\mu(B_x(r))\geq (nCr)^n\,.
	\end{equation*}
	In particular, for any $x\in X$, we have shown that
	\begin{equation*}
		\mu(B_x(r_0)))\geq\min\{\eta,\,(nCr_0)^n\}\,.
	\end{equation*}
	which proves the first property of the lemma.

	Next, dividing by $\omega_n r^n$ and taking the liminf when $r\rightarrow 0$ we obtain directly the second statement.
\end{proof}
It is rather easy to read on the isoperimetric profile whether a metric measure space is dominated by a space bilipschitz to a hyperbolic trumpet:
\begin{proposition}
Let $(X,d,\mu)$ be an infinite volume metric measure space.
Then it is dominated by a radial space radially bilipschitz to an $n$-dimensional hyperbolic trumpet
if and only if it has nonzero Cheeger constant and its isoperimetric profile $\Phi$ satisfies the asymptotic condition:
\begin{equation}
	\liminf_{t\rightarrow 0}\frac{\Phi(t)}{t^{1-\frac{1}{n}}}>0\,.
\end{equation}
\end{proposition}
\begin{proof}
	First assume that $(X,d,\mu)$ is dominated by a radial space $(\E^n,\sigma)$ radially $C$-bilipschitz to a
	hyperbolic trumpet $(\H^n_\alpha,\sigma=dr^2+\alpha^\frac{2}{n-1}\sinh^2(r)d\theta^2)$.
	Then a quick computation on the radial profile of the hyperbolic trumpet shows that
	\begin{equation*}
		\varphi^{\H^n_\alpha}(t)\underset{r\rightarrow 0}{\sim}\alpha^\frac{1}{n}n\omega_n^\frac{1}{n}t^{1-\frac{1}{n}}\,,\varphi^{\H^n_\alpha}(t)\geq t\,.
	\end{equation*}
	From the bilipschitz condition we deduce that
	\begin{equation*}
		\varphi^\sigma(t)\geq\frac{1}{C}\inf_{s\in[\frac{1}{C}t,Ct]}\varphi^{\H^n}_\alpha(s)
	\end{equation*}
	Since $\varphi^\sigma$ is a lower bound for the isoperimetric profile of $(X,d,\mu)$, we deduce the desired asymptotics.

	Conversely, if those two conditions on the isoperimetric profile are satisfied, there exists a continuous function $f$
	satisfying, for some $C>0$
	\begin{equation*}
		\lim_{t\rightarrow 0}\frac{f(t)}{t^{1-\frac{1}{n}}}>0\,,f(t)\geq C t,\, f(t)\leq\Phi(t)\,.
	\end{equation*}
	Such a function $f$ is always the radial profile of a space radially bilipschitz to a hyperbolic trumpet, as desired.
\end{proof}
\begin{remark}
	Note that we didn't discuss here whether the radial space dominating our metric measure space here is well-behaved or not.
	While it seems natural to ask for it, it is not needed if one wants to prove inequalities that remain true under bilipschitz transform, as is the case for the
	Moser--Trudinger inequality.
\end{remark}
\subsubsection{domination among $\rm{CD}(k,n)$-spaces}
Recall that the isoperimetric dimension and asymptotic isoperimetric ratio
have been defined in Definition~\ref{def:asymisolength}.
Among $\rm{CD}(k,n)$-spaces, we have the following characterization:
\begin{proposition}\label{prop:optimaldom}
	Let $(X,d,\mu)$ be a metric measure space satisfying a $\rm{CD}(k,n)$-condition. Let $m=n_X\leq n$ denote its isoperimetric dimension.
	Then it is dominated by a radial space bilipschitz to a hyperbolic trumpet if and only if its Cheeger constant is nonzero and $\calV(X)>0$.
	
	In that case, for any integer $l>m$, $(X,d,\mu)$ is dominated by an $l$-dimensional hyperbolic trumpet $\H^l_\theta$ for any cone angle $\theta>0$.

	If $m$ is an integer and its asymptotic isoperimetric ratio~$l_\infty(X)$ is nonnegative, then it is dominated by an $m$-dimensional hyperbolic trumpet
	$\H^m_{l_\infty}$.
\end{proposition}
\begin{proof}
	Let $(X,d,\mu)$ be such a space.
	First, since a space radially bilipschitz to a hyperbolic trumpet admits a nonzero cone singularity at zero, applying Lemma~\ref{lem:ballvolume}
	we obtain that if $(X,d,\mu)$ is dominated, $\calV(X)>0$.
	Also, the radial profile of a space $(\E^l,\sigma)$ $C$-bilipschitz to a hyperbolic trumpet $\H^n_\alpha$ satisfies
	\begin{equation*}
		\varphi^\sigma(t)\geq\frac{1}{C^2}\varphi^{\H^l_\alpha}(t)\geq \frac{1}{C^2}t\,,
	\end{equation*}
	Hence from $\Phi(t)\geq\varphi^\sigma(t)$ we deduce that the isoperimetric profile is surlinear, hence it admits a nonzero Cheeger constant.

	Conversely, let $(X,d,\mu)$ be an infinite volume $\rm{CD}(k,n)$-space with nonzero Cheeger constant $h>0$ and $\calV(X)>0$.
	By construction, for any $l>n_X$, $C>0$, there is $\eta>0$ such that
	\begin{equation*}
		0\leq t\leq\eta\,\Rightarrow\,\Phi(t)\geq Ct^{1-\frac{1}{l}}\,.
	\end{equation*}
	The Cheeger constant condition gives $h>0$ such that
	\begin{equation*}
		\Phi(t)\geq h\cdot t\,.
	\end{equation*}
	So there is a continuous function $f$ with the asymptotic $CT^{1-\frac{1}{l}}$ at zero, linear for $t$ large enough, such that
	\begin{equation*}
		f(t)\leq \Phi(t)\,.
	\end{equation*}
	Such a function is always the radial profile of a radial space radially bilipschitz to an $l$-dimensional hyperbolic trumpet with cone angle~$\frac{C^l}{l^l\omega_l^{l-1}}$ which dominates $(X,d,\mu)$, as desired. Since $C$ can be chosen arbitrarily, so can the cone angle of the hyperbolic trumpet.

	If $l=n_X$ is an integer and $l_\infty(X)>0$, the same construction yields a dominating radial space radially bilipschitw to an $l$-dimensional
	hyperbolic trumpet with cone angle $l_\infty(X)$, as claimed.
\end{proof}
\begin{remark}
	Note that for a noncompact finite volume $\rm{CD}(k,n)$-space, $\calV(X)=0$, hence it does not admit any domination.
	A compact $\rm{CD}(k,n)$-space always satisfies  $\calV(X)>0$, cf Lemma~\ref{lem:cpct}.
	Rajala~\cite{Raj12} showed that~$\rm{CD}(k,n)$-spaces supports the weak local Poincaré inequality,
	hence the Cheeger isoperimetric constant is nonzero. In particular, any compact~$\rm{CD}(k.n)$ has a nonzero Cheeger constant.
\end{remark}
\begin{corollary}\label{cor:domcompact}
	Let $(X,d,\mu)$ be a compact $\rm{CD}(k,n)$-space.
	Let $m=n_X$ denote its isoperimetric dimension
	Then for any $l>n_X$, $(X,d,\mu)$ is dominated by a radial space radially bilipschitz to an $l$-dimensional hyperbolic trumpet $\H^l_\theta$ for any cone angle $\theta>0$.
	If $n_X\in\N$ and its asymptotic isoperimetric ratio $l_\infty(X)$ is nonnegative, then it is dominated bt an $n_X$-dimensional hyperbolic trumpet
	with cone angle $l_\infty(X)$.
\end{corollary}
\begin{proof}
	We only have to show that $(X,d,\mu)$ has a nonzero Cheeger constant in order to apply the previous proposition.
	Let $(X,d,\mu)$ be such a space. By Rajala~\cite{Raj12}, it admits a $1-1$-local Poincaré inequality.
	Taking a radius large enough such that for any $x\in X$, $B_x(r)=X$, we obtain that there is $C>0$ such that, for any $u\in L^1(X)$:
	\begin{equation*}
		\int_X \big(u-\frac{1}{\mu(X)}\int_X ud\mu\big)d\mu\leq C\int_X\rm{lip}(u)d\mu\,.
	\end{equation*}
	Let $A\subset X$ a Borel set, such that $2\mu(A)\leq\mu(X)$.
	For any sequence $(f_i)$ converging to $\chi_A$ in $L^1(X)$, we can bound their variations:
	\begin{equation*}
		\liminf \int_X\rm{lip}(f_i)\geq C\lim\int\big(f_i-\frac{1}{\mu(X)}\int_X f_id\mu\big)d\mu\geq 2C\mu(A)(1-\frac{\mu(A)}{\mu(X)})\geq C\mu(A)\,.
	\end{equation*}
	In particular, Cheeger's isoperimetric constant $h$ satisfies $h\geq C>0$, as desired.
\end{proof}
\begin{proposition}
	Let $(X,d,\mu=\calH^n)$ be an infinite volume non-collapsed $\rm{RCD}(k,n)$-space.
	Assume $\calV(X)>0$ and it has nonzero Cheeger constant.
	Then it is dominated by a radial space radially bilipschitz to a hyperbolic trumpet $\H^n_\alpha$, with $\alpha=\theta_\infty$
	the asymptotic growth ratio of $(X,d,\mu)$.
\end{proposition}
\begin{proof}
	The only addition to the previous proposition is the interpretation of $\alpha$
	In this setup, the isoperimetric profile has been studied by Antonelli et al \cite{APP+24},
	they prove that the ratio $\frac{\Phi(t)}{t^{1-\frac{1}{n}}}$ admits a nonzero limit at zero $L$, which admits the following expression:
	\begin{equation*}
		L=n\bigg(\omega_n \underset{r\rightarrow 0}{\lim}\big(\underset{x\in M}{\inf}\frac{\mu(B_x(r))}{\omega_n r^n}\big)\bigg)^\frac{1}{n}=
		n(\omega_n\theta_\infty)^\frac{1}{n}\,.
	\end{equation*}
	From which we deduce that $\alpha=\theta_\infty$, as desired.
\end{proof}

\section{Application to the Moser--Trudinger inequality}
In this section, we consider the following Moser--Trudinger type inequality:
\begin{definition}
Let $(X,d,\mu)$ be a metric measure space of infinite volume, $\alpha>0$ and $n\geq 2$ an integer.
We say that $(X,d,\mu)$ satisfies the  Moser--Trudinger inequality $(MT)^n_\alpha$ if
there is a constant $C>0$ such that, for any $u\in W^{1,n}(X)$:
\begin{equation}
	\int \rm{lip}(u)^n\leq 1\,\Rightarrow \int \exp\big(\alpha u^\frac{n}{n-1}\big)
		-P_n(\alpha u^\frac{n}{n-1})d\mu\leq C\,.
\end{equation}
	where $P_n(t)=\sum_{j=0}^{n-2} \frac{t^j}{j!}$.
\end{definition}
\begin{remark}
	Since $X$ has infinite volume, in general there is no embedding $W^{1,n}(X)\rightarrow L^p(X)$
	for $p<n$, hence the necessity of the polynomial correction term.
	In this whole section, the notation $u^r$ is used to denote $|u|^r$, unless otherwise specified.
\end{remark}
\begin{remark}
	There exist $\rm{CD}(k,n)$-spaces for $n$ not an integer. It is natural to ask whether they would admit a finer version of the Moser--Trudinger inequality for non-integer
	$n$, but this is not in the scope of the paper here. 
\end{remark}
Satisfying such an inequality implies a nice geometry at small and large scale of the space.
That geometry will be encapsulated by two quantities, the Cheeger constant and the asymptotic growth ratio, see definitions~\ref{def:cheeger},\ref{def:agr}.
Denote also by $\lambda_p$ the spectral gap of the $p$-laplacian. Recall that by results
of Buser, Matei, $\lambda_p$ is controlled by the bound on the curvature and the Cheeger
constant for manifolds.
We will prove here the following result:
\begin{theorem}\label{thm:cara}
	Let $(X,d,\mu)$ be a locally compact geodesic $\rm{CD}(k,n)$-space with $\calV(X)>0$, and nonzero Cheeger constant.
	Let $n_X$ denote its isoperimetric constant and $l^{(m)}$ the $m$-th isoperimetric ratios.

	Then $(X,d,\mu)$ satisfies the Moser--Trudinger inequality $(MT)^m_\alpha$ for $m\geq n_X$ and $\alpha\leq m(l_\infty\omega_{m-1})^\frac{1}{m-1}$.
	Moreover, $(X,d,\mu)$ does not satisfy the Moser--Trudinger inequality $(MT)_\alpha^n$ for $\alpha> n(\theta_\infty\omega_{n-1})^\frac{1}{n-1}$.
\end{theorem}
\begin{remark}
	In the case of noncollapsed $\rm{RCD}(k,n)$-spaces, $n_X=n$ and $\theta_\infty=l_\infty$, so we have a complete answer to which Moser--Trudinger inequalities are satisfied.
\end{remark}
In the case of manifolds, this is actually a characterization of manifolds satisfying the Moser--Trudinger inequality among manifolds with Ricci curvature bounded from below.
\begin{theorem}\label{thm:caramanifold}
A complete Riemannian $n$-dimensional manifold $(M,g)$ with Ricci curvature bounded below
	satisfies the Moser--Trudinger inequality $(MT)^n_\alpha$ if and only if $\calV(M)>0$,
	its Cheeger constant is nonzero and $\alpha$ satisfies
	\begin{equation}
	\alpha\leq n(\theta_\infty\omega_{n-1})^\frac{1}{n-1}\,.
	\end{equation}
\end{theorem}
\subsection{Counterexamples to the Moser--Trudinger inequality}
\begin{lemma}
	Let $(X,d,\mu)$ be a doubling space with $\calV(X)>0$,
	i.e. such that
	\begin{equation}
		\inf_{x\in X}\mu(B_x(1))=0\,.
	\end{equation}
	Then for any $n\geq 2$ and $\alpha>0$, there is a sequence $u_n\in W^{1,n}(X)$ such that
	\begin{equation}
		\int\rm{lip}(u)^n\leq 1\quad\int_X \exp\big(\alpha u^\frac{n}{n-1}\big)-P_n(\alpha u^\frac{n}{n-1})
		\rightarrow+\infty\,.
	\end{equation}
\end{lemma}
\begin{proof}
	Since $(X,d,\mu)$ is assumed to be doubling, there is $C>0$
	such that, for any $x\in X$:
	\begin{equation*}
		\mu(B_x(2))\leq C\mu(B_x(1))\,.
	\end{equation*}
	Consider $r_m\rightarrow 0$, and consider a sequence $x_m\in X$ such that
	\begin{equation*}
		\mu(B_{x_m}(1))=r_m\,.
	\end{equation*}
	Consider the function $u_m\in W^{1,n}(X)$ defined by
	\begin{equation*}
		u_m(y)=\left\{\begin{array}{ll}
			T_m&\text{ if } d(x_m,y)\leq 1\\
			T_m(2-d(x,y))&\text{ if }d(x_m,y)\in[1,2]\\
			0&\text{otherwise}
		\end{array}\right.
	\end{equation*}
	for some $T_m>0$ to be specified later.

	Then we have the following estimates:
	\begin{equation*}
		\int\rm{lip}(u_m)^n\leq T_m^n\mu(B_{x_m}(2))-\mu(B_{x_m}(1))\leq
			T_m^n (C-1)r_m\,.
	\end{equation*}
	Choosing then
	\begin{equation*}
		T_m^n=\frac{1}{(C-1)r_m}\,
	\end{equation*}
	we obtain
	\begin{equation*}
		\|\rm{lip}(u_m)\|_n\leq 1\,,
	\end{equation*}
	while
	\begin{equation*}
		\int \exp\big(\alpha u_m^\frac{n}{n-1}\big)-P_n(\alpha u^\frac{n}{n-1})\geq
		r_m\bigg[\exp\big(\frac{\alpha}{r_m^\frac{1}{n-1}(C(k)-1)^\frac{1}{n-1}}\big)
			-P_n(\frac{\alpha}{r_m^\frac{1}{n-1}(C(k)-1)^\frac{1}{n-1}})\bigg]\,.
	\end{equation*}
	which is arbitrarily large when $r_m\rightarrow 0$, as desired.
\end{proof}
\begin{remark}
	From the Bishop--Gromov comparison theorem, any $\rm{CD}(k,n)$-space is doubling.
\end{remark}
\begin{lemma}
	Let $(M,g)$ be a complete Riemannian manifold with Ricci curvature bounded below.
	Assume its Cheeger constant $h$ is zero.
	Then it doesn't satisfy any $(MT)_\alpha$ for any $\alpha>0$.
\end{lemma}
\begin{proof}
	Let $(M,g)$ be with Ricci curvature bounded below, and zero Cheeger constant.
	The generalized Buser inequality, Theorem~\ref{thm:buser}, shows that for any $C>0$,
	there is a function $u\in W^{1,n}(M)$ satisfying
	\begin{equation*}
		\|\nabla u\|_n\leq 1 \text{ and } \|u\|_n\geq C\,.
	\end{equation*}
	In particular, $(M,g)$ does not satisfy the Moser--Trudinger inequality $(MT)^n_\alpha$
	for any $\alpha>0$.
\end{proof}
\begin{remark}
	The same proof shows that any metric measure space with zero $n$-th spectral gap
	cannot satisfy $(MT)^n_\alpha$ for any $\alpha>0$. What's missing in general
	is that a priori there could exist metric measure space with zero Cheeger
	constant, yet nonzero $n$-th spectral gap. De Ponti--Mondino~\cite{DM21} showed it doesn't exist in the case $n=2$ for $\rm{RCD}(k,2)$-spaces,
	relying on Heat flow arguments. It is unclear to the author if a generalized version of Buser's inequality exists in metric measure spaces.
\end{remark}
For $\rm{CD}(k,n)$-spaces, let's first give an upper bound on the validity of the Moser--Trudinger inequality:
\begin{lemma}\label{lem:noMTcone}
	Let $(X,d,\mu)$ be a $\rm{CD}(k,n)$-space. Denote
	by $\theta_\infty$ the asymptotic growth of balls
	\begin{equation}
		\theta_\infty=\lim_{r\rightarrow 0}\inf_{x\in M}\frac{\Vol(B_x(r))}{\omega_nr^n}\,.
	\end{equation}
	Then for any $\alpha>n(\theta_\infty\omega_{n-1})^\frac{1}{n-1}$,
	there is a sequence $(u_m)\in W^{1,p}(X)$ such that
	\begin{equation}
		\|\nabla u_m\|_n\leq 1\text{ and }\int_X e^{\alpha u_m^\frac{n}{n-1}}
		-P_n(\alpha u_m^\frac{n}{n-1})\rightarrow+\infty\,.
	\end{equation}
\end{lemma}
\begin{proof}
	It mainly comes from the Bishop--Gromov comparison theorem. 
	Denote by $b(k,r)$ the volume of a radius $r$ ball in the constant curvature $k$ $n$-dimensional space, and $s(k,r)=\partial_r b(k,r)$ the area
	of spheres in that same space.
	Fix $\theta>\theta_\infty$, $r>0$ and $x_r\in X$ such that
	\begin{equation*}
		\mu(B_x(r))=\theta b(k,r)\,.
	\end{equation*}
	From the Bishop--Gromov comparison theorem, we deduce that for any $s>r$
	\begin{equation*}
		\mu(B_x(s))\leq\theta b(k,s)
	\end{equation*}
	Hence
	\begin{equation*}
		\frac{\mu(B_x(s))-\mu(B_x(r))}{s-r}\leq\frac{\theta b(k,s)-\theta b(k,r)}{s-r}
	\end{equation*}
	from which we deduce
	\begin{equation*}
		\rm{Per}(S_x(r))\leq \theta s(k,r)\,.
	\end{equation*}
	And again, from the Bishop--Gromov comparison theorem, for any $s\geq r$,
	\begin{equation*}
		\rm{Per}(S_x(s))\leq\theta s(k,s)\,.
	\end{equation*}
	To use the counterexamples used by Moser in~\cite{Mos71}, we work at a scale
	so that the curvature doesn't play too much of a role. Let $R>0$ and $\eta>0$
	be so that
	\begin{equation*}
		\forall s\in (0,R),\,1-\eta\leq\frac{b(k,s)}{b(0,s)}\leq 1+\eta\,.
	\end{equation*}
	Let $u:X\rightarrow \R$ be radially symmetric around $x$.
	Then one has, using $r$ as the radial geodesic coordinate around $x$,
	and introduce $t$ the variable such that $e^{-t}=\frac{r^n}{R^n}$:
	\begin{flalign*}
		\int_X \rm{lip}(u)^nd\mu&=\int_{\R^+}\rm{Per}(S(x,r))|\partial_r u|^n dr\\
			&=\int_{\R_+}n^{n-1}\rm{Per}(S(x,e^{-\frac{t}{n}}))|\partial_t u|^n
			e^{\frac{n-1}{n}t}dt
	\end{flalign*}
	Consider now the function $u$ radial around $x_r$ given by the following, and denote
	$t_0=-n\ln(\frac{r}{R})$.
	\begin{equation*}
		u(t)=\left\{\begin{array}{ll}
			C t &\text{ if } 0\leq t\leq t_0\\
			C t_0 &\text{ if }  t_0\leq t
		\end{array}\right.
	\end{equation*}
		whose support is $B(x_r,R)$. Then we can compute the norm of its gradient
	\begin{flalign*}
		\|\rm{lip}(u)\|_n^n&=\int_0^{t_0}n^{n-1}\rm{Per}(S(x_r,e^{-\frac{t}{n}}))C^n
			e^{\frac{n-1}{n}t}dt\\
			&\leq
			C^n n^{n-1} \theta(1+\eta)\omega_{n-1}t_0
	\end{flalign*}
	While we have the lower bound
	\begin{flalign*}
		I(u)&=\int_X\exp\big(\alpha u^\frac{n}{n-1}\big)-P_n(\alpha u^\frac{n}{n-1})d\mu\geq
		\mu(B(x,r))\big[\exp(\alpha C^\frac{n}{n-1} t_0^\frac{n}{n-1})
			-P_n(\alpha C^\frac{n}{n-1} t_0^\frac{n}{n-1})\big]\\
		&\geq\theta(1-\eta)\omega_{n-1}\frac{r^n}{n}\big[\exp(\alpha C^\frac{n}{n-1} 
			t_0^\frac{n}{n-1})-P_n(\alpha C^\frac{n}{n-1} t_0^\frac{n}{n-1})\big]\\
		&\geq
		\theta(1-\eta)\omega_{n-1}\frac{r^n}{n}\big[\exp(\alpha C^\frac{n}{n-1}
			n^\frac{n}{n-1}|\ln(\frac{r}{R})|^\frac{n}{n-1})
			-P_n(\alpha C^\frac{n}{n-1}
				n^\frac{n}{n-1}|\ln(\frac{r}{R})|^\frac{n}{n-1})\big]
	\end{flalign*}
	Choosing for $C$ the value
	\begin{equation*}
		C=\frac{1}{n\theta^\frac{1}{n}(1+\eta)^\frac{1}{n}
			\omega_{n-1}^\frac{1}{n}|\ln(\frac{r}{R})|^\frac{1}{n}}
	\end{equation*}
	we obtain the controls
	\begin{flalign*}
		\|\rm{lip}(u)\|_n&\leq 1\\
		I(u)&\geq \theta(1-\eta)\omega_{n-1}\frac{r^n}{n}\big[\exp(
			\frac{\alpha}{(\theta\omega_{n-1}(1+\eta))^\frac{1}{n-1}}
			|\ln(\frac{r}{R})|)
			-P_n(\frac{\alpha}{(\theta\omega_{n-1}(1+\eta))^\frac{1}{n-1}}
			|\ln(\frac{r}{R})|)\big]\,.
	\end{flalign*}
	When $r\rightarrow 0$, the quantity $I(u)$ diverges as soon as
	\begin{equation*}
		n<\frac{\alpha}{(\theta\omega_{n-1}(1+\eta))^\frac{1}{n-1}}\,.
	\end{equation*}
	Since $\theta$ can be chosen arbitrarily close to $\theta_\infty$ and $\eta$ arbitrarily close
	to $1$, this shows that for any $\alpha$ satisfying
	\begin{equation*}
		\alpha>n(\theta_\infty\omega_{n-1})^\frac{1}{n-1},
	\end{equation*}
	there is a sequence $(u_m)\in W^{1,n}(X)$ such that
	\begin{equation*}
		\|\rm{lip}(u_m)\|_n\leq 1\text{ and }\int_X \exp(\alpha u_m^\frac{n}{n-1})-
		P_n(\alpha u^\frac{n}{n-1})d\mu\rightarrow+\infty\,,
	\end{equation*}
	as claimed.
\end{proof}
\subsection{Positive results with regard to the Moser--Trudinger inequality}
We will deduce all our results from the study of the Moser--Trudinger inequality in the hyperbolic space,
which was covered in Mancini--Sandeep~\cite{MS10} for the plane,
and Mancini--Sandeep--Tintarev~\cite{MST13} in the higher-dimensional case.
\begin{theorem}[Mancini--Sandeep, Mancini--Sandeep--Tintarev]
	Consider $\H^n$ the $n$-dimensional hyperbolic space.
	There is a constant $C(n)>0$ such that, for any $u\in W^{1,n}(\H^n)$, the
	following implication is true:
	\begin{equation}
		\|\nabla u\|_n\leq 1\,\Rightarrow\,
		\int_{\H^n}\exp(\alpha_n u^\frac{n}{n-1})-P_n(\alpha u^\frac{n-1}{n})\leq C(n)\,.
	\end{equation}
	where $\alpha=\alpha_n=n\omega_{n-1}^\frac{1}{n-1}$.
\end{theorem}
A first corollary concerns the hyperbolic trumpets with cone singularity at zero.
\begin{corollary}\label{cor:MTtrumpet}
	Let $\beta\in (0,1)$ and 
	let $(\E^n, dr^2+\beta^\frac{2}{n-1}\sinh(r)^2d\theta^2)$ be the hyperbolic trumpet of cone angle $\beta$.
	Then for any $\alpha\leq n(\beta\omega_{n-1})^\frac{1}{n-1}$,
	the hyperbolic trumpet satisfies $(MT)^n_\alpha$
\end{corollary}
\begin{proof}
Fix $\alpha\leq n\beta\omega_{n-1}^\frac{1}{n-1}$.
First, since $\beta<1$ the hyperbolic trumpet of cone angle $\beta$ is well-behaved by Lemma~\ref{lem:wellbehaved}.
Hence the rearrangement map $u\mapsto \hat u$ on itself satisfies:
\begin{equation*}
	\|\nabla \hat u\|_n\leq\|\nabla u\|_n\,\text{and}\,
	\int e^{\alpha u^\frac{n}{n-1}}-P_n(\alpha u^\frac{n}{n-1})=
	\int e^{\alpha \hat u^\frac{n}{n-1}}-P_n(\alpha \hat u^\frac{n}{n-1})\,.
\end{equation*}
In other words, it is enough to prove the Moser--Trudinger inequality on radial functions
on the hyperbolic trumpet.
Now consider the identity map
	\begin{equation*}
		F:\left\{\begin{array}{ll}
			W^{1,n}_{rad}(\H^n_\beta)&\rightarrow W^{1,n}(\H^n)\\
			u&\mapsto u\end{array}\right.
	\end{equation*}
A quick computation shows that
		\begin{equation*}
			\|\nabla u\|_n=\beta^\frac{1}{n}\|\nabla F(u)\|_n\text{ and }
			\int_{\H^n_\beta} e^{\alpha u^\frac{n}{n-1}}-P_n(\alpha u^\frac{n}{n-1})
			=\beta\int_{\H^n} e^{\alpha F(u)^\frac{n}{n-1}}-P_n(\alpha F(u)^\frac{n}{n-1})
		\end{equation*}
	Assume that $u\in W^{1,n}_{rad}(\H^n_\beta)$ satisfies $\|\nabla u\|_n\leq 1$.
	Denoting $v=\beta^\frac{1}{n}F(u)\in W^{1,n}(\H^n)$, it satisfies
	\begin{equation*}
		\|\nabla v\|_n\leq 1\,.
	\end{equation*}
	Also, by assumption
	\begin{equation*}
		\alpha \beta^\frac{-1}{n-1}\leq n\omega_{n-1}^\frac{1}{n-1}=\alpha_n\,. 
	\end{equation*}
	By the Mancini--Sandeep--Tintarev Theorem, there is a constant $C(n)$
	such that
	\begin{equation*}
		\int_{\H^n} e^{\alpha\beta^\frac{-1}{n-1}v^\frac{n}{n-1}}
			-P_n(\alpha\beta^\frac{-1}{n-1}v)\leq C(n)
	\end{equation*}
	i.e.
	\begin{equation*}
		\int_{\H^n_\beta} e^{\alpha u^\frac{n}{n-1}}-P_n(\alpha u^\frac{n}{n-1})
		\leq C(n)\,.
	\end{equation*}
	So the hyperbolic trumpet $\H^n_\beta$ satisfies $(MT)^n_\alpha$ for
	$\alpha\leq n(\beta\omega_{n-1})^\frac{1}{n-1}$, as claimed.
\end{proof}
Satisfying the Moser--Trudinger inequality is invariant under Bilipschitz conformal
transforms, as the following lemma states:
\begin{lemma}\label{lem:MTbilip}
	Let $(M_1,g_1)$ and $(M_2,g_2)$ be complete Riemannian manifolds of dimension $n$,
	and assume that there exists a bilipschitz conformal homeomorphism $f:M_1\rightarrow M_2$.
	Then for $\alpha>0$, $(M_1,g_1)$ satisfies $(MT)^n_\alpha$
	if and only if $(M_2,g_2)$ satisfies $(MT)^n_\alpha$.
\end{lemma}
\begin{proof}
	Since the statement is symmetric, it is enough to assume that $(M_1,g_1)$
	satisfies $(MT)_\alpha$ and to prove that $(M_2,g_2)$ satisfies it too.
	Let $C>0$ be such that, for any $u\in W^{1,n}(M_1)$:
	\begin{equation*}
		\|\nabla u\|_n^n\leq 1\,\Rightarrow\,
		\int e^{\alpha u^\frac{n}{n-1}}-P_n(\alpha u^\frac{n}{n-1})\leq C\,.
	\end{equation*}
	Let $v\in W^{1,n}(M_2)$ such that $\|\nabla v\|_n\leq 1$.
	Since $f$ is conformal, we have $\|\nabla(v\circ f)\|_n=\|\nabla v\|_n\leq 1$.
	And since it is bilipschitz conformal, we can control the integral
	\begin{equation*}
		\int_{M_2} e^{\alpha v^\frac{n}{n-1}}-P_n(\alpha v^\frac{n}{n-1})
		\leq \rm{Lip}(f)\int_{M_1}
			e^{\alpha (v\circ f)^\frac{n}{n-1}}-P_n(\alpha(v\circ f)^\frac{n}{n-1})
			\leq \rm{Lip}(f) C\,.
	\end{equation*}
	Hence $(M_2,g_2)$ satisfies $(MT)_\alpha$, as claimed.
\end{proof}
Finally, we are able to conclude with the proofs of Theorems~\ref{thm:cara},\ref{thm:caramanifold}.
\begin{proof}[Proof of Theorems~\ref{thm:cara},\ref{thm:caramanifold}]
	Let $(M,g)$ be an infinite volume
	complete riemannian $n$-manifold with Ricci curvature bounded below by $k\in\R$.
	We have already seen that if there is no lower bound on the volume of balls, then $(MT)^n_\alpha$ does not hold
	for any $\alpha>0$. Similarly, if its Cheeger constant is zero, then it cannot
	satisfy any Moser--Trudinger type inequality.

	Let now $(X,d,\mu)$ be a locally compact, geodesic $\rm{CD}(k,n)$-space of infinite volume.
	Assume it has nonzero Cheeger constant and $\calV(X)>0$.
	Let $m\geq n_X$ and $l\leq l_\infty^{(m)}(X)$, $\alpha=m(l_\infty^{(m)}\omega_{m-1})^\frac{1}{m-1}$.
	We know from Proposition~\ref{prop:optimaldom} that $(X,d,\mu)$ is dominated
	by a radial ball bilipschitz conformal to an $m$-dimensional hyperbolic trumpet of cone angle $l$.
	By Lemma~\ref{lem:MTbilip} and Corollary~\ref{cor:MTtrumpet},
	that dominating space satisfies the Moser--Trudinger inequality $(MT)^m_\alpha$
	Since this space dominates $(X,d,\mu)$, the same is true for the space $(X,d,\mu)$,
	as desired. Note that this is the only point where we used the asumption that $(X,d,\mu)$ is locally compact and geodesic.
	
	For $\alpha>n(\theta_\infty\omega_{n-1})^\frac{1}{n-1}$, it does not satisfy $(MT)^n_\alpha$
	from Lemma~\ref{lem:noMTcone}, hence the characterization is proven.
\end{proof}
\subsection{Dependency of the Moser--Trudinger constant}
The method used here allows to give a uniform Moser--Trudinger inequality for families of
locally compact geodesic metric measure spaces $(X,d,\mu)$ dominated by a common space bilipschitzly conformal to a hyperbolic trumpet.
In particular it allows the specific case:
\begin{theorem}
	Let $a,b\in\R$, and $\delta>0$ and $h_0>0$
	There is a constant $C=C(n,a,b,\delta,h_0)$ such that
	for any $n$-dimensional infinite volume manifold satisfying
	\begin{equation}
		a\leq\rm{Sect}(M)\leq b,\,\rm{inj}(M)\geq \delta,\, h(M)\geq h_0\,,
	\end{equation}
	then for $\alpha\leq n\omega_{n-1}^\frac{1}{n-1}$ and $u\in W^{1,n}(M)$ one has
	\begin{equation}
		\|\nabla u\|_n\leq 1\,\Rightarrow
		\int_M \exp\big(\alpha u^\frac{n}{n-1}\big)-P_n(\alpha u^\frac{n}{n-1})\leq C\,.
	\end{equation}
\end{theorem}
\begin{proof}
	The controls on the sectional curvature are enough, by the Bishop--Gunther and Bishop--Gromov comparison theorems, to control the isoperimetric profile of the manifolds uniformly. Hence at fixed $a,b,\delta,h_0$, there is a common radial ball dominating all
	the manifolds satisfying the assumptions, so there is a uniform constant
	valid for all the inequalities $(MT)^n_\alpha$ involved.
\end{proof}

\section{Double rearrangement for compact manifolds}
In this section we consider compact $\rm{CD}(k,n)$-spaces,
and introduce a double rearrangement
specifically adapted to functions with zero average.
As a consequence, we are able to show that
all compact manifolds admit a Moser--Trudinger type
inequality, and we discuss its sharpness and
dependency in the geometry of the manifold.
\subsection{Double rearrangement for $\rm{CD}(k,n)$-spaces}
We already introduced a notion of domination adapted to infinite volume manifolds.
Here we use this notion in order to introduce a double rearrangement for functions on a finite volume $\rm{CD}(k,n)$-space.
In all that follows $(X,d,\mu)$ denotes a finite volume, locally compact, geodesic metric measure space.
We will make use of a particular value called the median of $u$:
\begin{definition}
	Let $(X,d,\mu)$ be a finite volume metric measure space.
	Let $u:X\rightarrow\R$ be integrable.
	We call median of $u$, and denote by $c(u)$ or $c$ when no confusion is possible,
	any value $c\in\R$ satisfying:
	\begin{equation}
		\mu(\{u<c\})\leq\frac{\mu(X)}{2}\,\text{and}\,\mu(\{u>c\})\leq\frac{\mu(X)}{2}\,.
	\end{equation}
\end{definition}
\begin{remark}
	Both these inequalities are equalities as soon as $\mu(\{u=c\})=0$.
	While $c$ may not be uniquely defined, the set of medians of $u$ is always an interval,
	so one can choose any particular point in that interval, as it won't play a role
	in the forthcoming discussion.
\end{remark}
For $u\in L^p(X)$, its median and average are closely related by the Markov's inequality
\begin{lemma}\label{lem:diffavgmed}
	Let $u\in L^p(X)$. Denote by $\overline u$ its average. Then any median $c$ of $u$ satisfies:
	\begin{equation}
		|c-\overline u|\leq(\frac{2}{\mu(X)})^\frac{1}{p}\|u-\overline u\|_p\,.
	\end{equation}
\end{lemma}
\begin{proof}
	From Markov's inequality, for any $t>0$, we have
	\begin{equation*}
		\mu(\{|u-\overline u|\geq t\})\leq\frac{\|u-\overline u\|_p^p}{t^p}\,.
	\end{equation*}
	Fix
	\begin{equation*}
		t=(\frac{2}{\mu(X)})^\frac{1}{p}\|u-\overline u\|_p\,.
	\end{equation*}
	Then from Markov's inequality, we obtain, for any $\eps>0$:
	\begin{equation*}
		\mu(\{u<t+\overline u+\eps\})>\frac{\mu(X)}{2}\text{ and }
		\mu(\{u>\overline u-t-\eps\})>\frac{\mu(X)}{2}\,.
	\end{equation*}
	From which we deduce that any median $c$ satisfies
	\begin{equation*}
		|c-\overline u|\leq t+\eps\,,
	\end{equation*}
	This being true for any choice of $\eps>0$, we have proven the desired result.
\end{proof}
We will also need a lower estimate on the gradient of some functions, relying on
the spectral gap of the $n$-Laplacian of the space.
\begin{lemma}\label{lem:bettergradient}
	Let $\beta<1$, and consider $\H^n_\beta$ the hyperbolic trumpet with cone angle $\beta$.
	Let $\Omega\subset\H^n_\beta$ a bounded domain.
	Let $c>0$ and $u:\H^n_\beta\rightarrow\R$ be a function such that $u-c$ is supported in $\Omega$
	and satisfies
	\begin{equation}
		u\leq c\text{ and }\int_\Omega u\leq -2c\sigma(\Omega)\,.
	\end{equation}
	Then it satisfies
	\begin{equation}
		\|\nabla u\|^n_n\geq\beta^\frac{1}{n-1}(3(1-\frac{1}{n}))^nc^n\sigma(\Omega)\,.
	\end{equation}
\end{lemma}
\begin{proof}
	First apply the Hölder inequality to obtain, for any $1\leq k\leq n$:
	\begin{equation*}
		|\int_\Omega u|\leq\|u|_\Omega\|_k^k\sigma(\Omega)^\frac{k-1}{k}\,.
	\end{equation*}
	Hence we deduce
	\begin{equation*}
		2^k c^k\sigma(\Omega)\leq\|u|_\Omega\|_k^k
	\end{equation*}
	We estimate now the norm of $u-c$:
	\begin{equation*}
		\|u-c\|_n^n=\sum_{k=0}^n\binom{n}{k}\|u|_\Omega\|_k^k c^{n-k}
			\geq\sum_{k=0}^n\binom{n}{k}2^kc^n\sigma(\Omega)\geq 3^nc^n\sigma(\Omega)\,.
	\end{equation*}
	Note that we used the assumption $u\leq c$ here.
	Finally, since $u-c$ is compactly supported, we obtain from Buser's inequality
	\begin{equation*}
		\|\nabla u\|_n^n\geq\mu_n(\H^n_\beta)3^n c^n\sigma(\Omega)
	\end{equation*}
	It remains to estimate $\mu_n(\H^n_\beta)$.
	From the variational expression
	\begin{equation*}
		\mu_n(M)=\inf\big\{\frac{\int_M\|\nabla f\|^n}{\int_M f^n}, f\in W^{1,n}(M)\big\}\,,
	\end{equation*}
	we derive that
	\begin{equation*}
		\mu_n(\H^n_\beta)\geq \beta^\frac{1}{n-1}\mu_n(\H^n)\geq \beta^\frac{1}{n-1}\frac{(n-1)^n}{n^n}\,.
	\end{equation*}
	Where we used the estimate for $\mu_n(\H^n)$ from Lima--Montenegro--Santos~\cite{LMS10}.
	We have proven that
	\begin{equation*}
		\|\nabla u\|^n_n\geq\beta^\frac{1}{n-1}(3(1-\frac{1}{n}))^nc^n\sigma(\Omega)\,,
	\end{equation*}
	as claimed.
\end{proof}
For $u\in W^{1,p}(X)$, and $c\in\R$ such that $\mu(\{u>c\})\leq\frac{\mu(X)}{2}$ and
$\mu(\{u<c\})\leq\frac{\mu(X)}{2}$, we consider $u_+$ and $u_-$ the rearrangements
of $(u-c)\one_{u>c}$ and of its negative part, in the following way
\begin{definition}
	Let $X$ be a compact $\rm{CD}(k,n)$-space dominated by a radial space $(\E^m,\sigma)$,
	and $p>1$.
	For $u\in W^{1,p}(X)$ and $c\in\R$ its median, we consider $u_+,u_-\in L^p(X)$
	the radially symmetric decreasing functions satisfying:
	\begin{equation}
		\left\{\begin{array}{ll}
			\sigma(\{u_+> t\})=\mu(\{u-c>t\})&\forall t>0\,\inf u_+=0\\
			\sigma(\{u_->t\})=\mu(\{c-u>t\})&\forall t>0\,\inf u_-=0\,.
		\end{array}\right.
	\end{equation}
\end{definition}
The worth of this rearrangement is in the following properties:
\begin{proposition}\label{PScompact}
	For $u\in W^{1,p}(X)$ with median $c$, consider $u_-,u_+$ its rearrangements on $(\E^m,\sigma)$.
	Let $F:\R\rightarrow\R$ be a measurable function, such that $F(c)=0$. Then we have
	\begin{equation}
		\int_X F(u)d\mu=\int_{\E^m}F(u_++c)d\sigma+\int_{\E^m} F(c-u_-)d\sigma\,.
	\end{equation}
	As with the following Polya--Szegö type inequality:
	\begin{equation}
		\|\nabla u_+\|^p_p+\|\nabla u_-\|_p^p\leq\|\rm{lip}(u)\|_p^p\,.
	\end{equation}
\end{proposition}
\begin{proof}
	Let $\chi(t)=\max(t,0)$.
	Then decompose $u$ in the following way:
	\begin{equation*}
		u=\chi(u-c)+c-\chi(c-u)\,.
	\end{equation*}
	From that decomposition, we recognize
	that $u_+$ is the rearrangement of $\chi(u-c,0)$ while $u_-$ is the rearrangement
	of~$\chi(c-u)$.
	Since $F(c)=0$, we can write the integral:
	\begin{flalign*}
		\int_X F(u)&=\int_X F(\chi(u-c)+c)d\mu+\int_X F(c-\chi(c-u))d\mu\\
			&=\int_{\E^m} F(u_++c)d\sigma+\int_{\E^m} F(c-u_-)d\sigma\,.
	\end{flalign*}
	From the same decomposition. we can write the gradient as a sum:
	\begin{equation*}
		\rm{lip}(u)=\rm{lip}(\chi(u-c))+\rm{lip}(\chi(c-u))\,.
	\end{equation*}
	Applying the Polya-Szegö inequality, valid since we have the required domination property,
	we obtain
	\begin{equation*}
		\|\nabla u_+\|_p^p+\|\nabla u_-\|_p^p\leq
			\int_X \rm{lip}(\chi(u-c))^p+\rm{lip}(\chi(c-u))^pd\mu\leq\|\rm{lip}(u)\|_p^p\,,
	\end{equation*}
	as claimed.
\end{proof}
\subsection{Moser--Trudinger inequality for compact metric measure spaces}
In the case of compact $\rm{CD}(k,n)$-spaces, we are able to prove the
following Moser--Trudinger type inequality:
\begin{theorem}\label{thm:MTcompact}
	Let $(X,d,\mu)$ be a compact $\rm{CD}(k,n)$-space whose volume is larger than $V_0>0$
	Let $n_X$ be its isoperimetric dimension, and $l_\infty$ its asymptotic isoperimetric ratio, defined by:
	\begin{equation}
		l_\infty=\liminf_{t\rightarrow 0}\frac{\Phi(t)^m}{m^m\omega_m T^{m-1}}\,.
	\end{equation}
	and by $\varphi$ the radial profile of a radial space dominating $X$, radially bilipschitzly
	conformal to $\H^m_\beta$, with $\omega_{m-1}\beta^{m-1}=l_\infty$.
	Then for any $m>n_X$ and $\alpha>0$ or for $m=n_X$ and $\alpha\leq m(l_\infty)^\frac{1}{m-1}$,
	there is a constant $C(m,\alpha,V_0,\varphi)$ such that, for any $u\in W^{1,m}_0(X)$,
	we have
	\begin{equation*}
		\|\rm{lip}(u)\|_m\leq 1\,\Rightarrow\,\int_X e^{\alpha u^\frac{m}{m-1}}
		-P_m(\alpha u^\frac{m}{m-1})d\mu \leq C(m,\alpha,V_0,\varphi)\,.
	\end{equation*}
\end{theorem}
We will deduce this Theorem from the noncompact case, via double rearrangement.

The main idea is that compact $\rm{CD}(k,n)$-spaces are dominated by radial spaces, which are bilipschitz to hyperbolic trumpets, from Corollary~\ref{cor:domcompact}.

The first step of the proof is to replace $W^{1,m}(X)$ with the space of zero median
functions, more adapted to double rearrangements.
\begin{proposition}\label{prop:Mtmed}
	Let $(X,d,\mu)$ be a compact $\rm{CD}(k,n)$-space.
	For $u\in W^{1,m}(X)$, denote by $c$ a median of $u$.
	Then for any $\alpha\leq m(l_\infty^{(m)}\omega_{m-1})^\frac{1}{m-1}$,
	there is a constant $C(m,\alpha)$ such that
	\begin{equation}
		\|\rm{lip}(u)\|_m\leq 1\,\Rightarrow\,
		\int_X e^{\alpha (u-c)^\frac{m}{m-1}}-P_m(\alpha (u-c)^\frac{m}{m-1})\leq C(m,\alpha)\,.
	\end{equation}
	The constant $C(m,\alpha)$ can be replaced by a constant depending only on
	the isoperimetric profile of $(X,d,\mu)$.
\end{proposition}
\begin{proof}
	From Corollary~\ref{cor:domcompact} we know that $(X,d,\mu)$ is dominated by a radial space $(\E^m,\sigma)$ radially bilipschitz
	conformal to a hyperbolic trumpet of cone angle $l\leq l^{(m)}_\infty$.
	Let $\alpha\leq m(l^{(m)}_\infty\omega_{m-1})^\frac{1}{m-1}$.
	We know that there is a constant $C>0$ such that, for $u\in W^{1,m}(\E^m)$,
	\begin{equation*}
		\|\nabla u\|_m\leq 1\,\Rightarrow
		\int_{\E^m} e^{\alpha u^\frac{m}{m-1}}-P_m(\alpha u^\frac{m}{m-1})d\sigma\leq C\,.
	\end{equation*}
	For $u\in W^{1,m}(X)$ with $\|\rm{lip}(u_)\|_m\leq 1$ and median $c$,
	consider its double rearrangement $u_-$,$u_+$. From Proposition~\ref{PScompact},
	we know that $u_-,u_+\in W^{1,m}(\E^m)$, with $\|\nabla u_-\|_m^m+\|\nabla u_+\|_m^m\leq 1$.
	Applying the Moser--Trudinger inequality to both of these, we obtain
	\begin{equation*}
		\int_X e^{\alpha (u-c)^\frac{m}{m-1}}-P_m(\alpha (u-c)^\frac{m}{m-1})d\mu\leq 2C\,.
	\end{equation*}
	as wished.
\end{proof}
To replace the median by the average in the Moser--Trudinger inequality, the proof is a bit
more involved. Let's introduce the notations first
\begin{notation}
	We will make use of the following notations:
	\begin{enumerate}
		\item $(X,d,\mu)$ is a finite volume $\rm{CD}(k,n)$-space dominated by a radial space $(\E^m,\sigma)$ with radial profile $\varphi$,
			radially bilipschitz conformal to a hyperbolic trumpet $\H^m_\beta$.
			The volume of $X$ is assumed to be larger than $V_0>0$.
		\item $u$ is a function on $X$ with zero average, and positive median $c>0$,
			such that $\|\rm{lip}(u)\|_m\leq 1$.
		\item $u_-,u_+$ are the rearrangements of $u$ on $(\E^m,\sigma)$. We denote by $\Omega\subset\E^m$ the centered ball of volume half the volume of $X$.
		\item $F$ denotes the following function:
			\begin{equation}
				F(t)=e^{|t|}-\sum_{j=0}^{m-2}\frac{|t|^j}{j!}\,.
			\end{equation}
		\item $l_\infty^{(m)}$ is the $m$-th asymptotic isoperimetric ratio of $X$,
			$\alpha$ is a constant such that $\alpha\leq m(l_\infty^{(m)}\omega_{m-1})^\frac{1}{m-1}$.
	\end{enumerate}
\end{notation}
\textbf{Step 1:The Moser Trudinger inequality for $u_-$}
\begin{proposition}\label{prop:step1}
	With the notations introduced, there is a constant $C_5(\varphi,\alpha,h,V_0)$
	such that $u_-$ satisfies
	\begin{equation}
		\int_\Omega F(\alpha(c-u_-)^\frac{m}{m-1})d\sigma\leq C_5(\varphi,\alpha,h,V_0)\,.
	\end{equation}
	Moreover, the constant can be chosen so that $C_5$ is uniform for all
	volumes larger than $V_0>0$.
\end{proposition}
\begin{proof}
	First, from the Moser--Trudinger inequality for zero-median functions,
	we get $C_1>0$ such that
	\begin{equation*}
		\int_X F(\alpha(u-c)^\frac{m}{m-1})d\mu\leq C_1(\varphi,\alpha)\,.
	\end{equation*}
	From the properties of the double rearrangements, as $G=F(\alpha(\cdot-c)^\frac{m}{m-1})$
	vanishes at $c$, we have
	\begin{equation*}
		\int_X F(\alpha(u-c)^\frac{m}{m-1})d\mu=\int_{\E^m} F(\alpha u_+^\frac{m}{m-1})d\sigma+
					\int_{\E^m} F(\alpha u_-^\frac{m}{m-1})d\sigma\,.
	\end{equation*}
	So in particular,
	\begin{equation*}
		\int_{\E^m} F(\alpha u_-^\frac{m}{m-1})d\sigma\leq C_1\,.
	\end{equation*}
	By construction, $u_-$ is supported in $\Omega$  and $0\leq c-u_-\leq c$.
	In particular, when $u_-\geq\frac{c}{2}$ we have
	$|c-u_-|\leq|u_-|$, hence we get the following bounds:
	\begin{equation*}
		\int_\Omega F(\alpha(u_--c)^\frac{m}{m-1})d\sigma
		\leq \int_{\E^m} F(\alpha(u_-)^\frac{m}{m-1})d\sigma
		+V_0\underset{[\frac{c}{2},c]}{\sup} F(\alpha(\cdot)^\frac{m}{m-1})\,.
	\end{equation*}
	From the Markov inequality, we know that
	\begin{equation*}
		|c|\leq\frac{C_2}{V_0}\|u\|_2\leq\frac{C_3(h)}{V_0}\,.
	\end{equation*}
	Where we used that a control on $h$ implies a control on the $2$-spectral gap, cf De Ponti--Mondino~\cite{DM21}.
	A study of the variations of $F$ then shows that we can estimate brutally the supremum:
	\begin{equation*}
		V_0\underset{[\frac{c}{2},c]}{\sup} F(\alpha(\cdot)^\frac{m}{m-1})\leq
		\sum_{j=m-1}^\infty\frac{\alpha^j}{j!} C_3^\frac{jm}{m-1}V_0^{1-\frac{jm}{m-1}}\leq C_4(\alpha,h,V_0)\,.
	\end{equation*}
	Hence there is $C_5(\varphi,\alpha,h,V_0)$ such that
	\begin{equation*}
		\int_\Omega F(\alpha(c-u_-)^\frac{m}{m-1})d\sigma\leq
			C_5(\varphi,\alpha,h,V_0)\,,
	\end{equation*}
	as claimed. Note that the constant appearing gets arbitrarily large only when $V_0$ is close to zero,
	but not when $V_0$ is large, hence the comment on the dependency of $C_4$.
\end{proof}
\textbf{Step 2: get an improved bound on $\nabla u_+$}
\begin{proposition}\label{prop:step2}
	With the notations introduced, there is a constant $C_6=C_6(\varphi)$
	such that $u_+$ satisfies
	\begin{equation*}
		\|\nabla u_+\|_m^m\leq 1-C_6c^m\sigma(\Omega)\,.
	\end{equation*}
\end{proposition}
\begin{proof}
	Let $\beta=\beta(\varphi)$ be such that the radial space dominating
	$(X,d,\mu)$ is bilipschitzly conformal to the trumpet $\H^m_\beta$.

	Let $\Psi:\H^m_\beta\rightarrow\E^m$ be a radially biLipschitz conformal map.
	The function $c-u_-$ satisfies:
	\begin{equation*}
		\rm{Supp}(c-u_-)\subset \Omega\text{ and } w\leq c
	\end{equation*}
	Also, by construction
	\begin{equation*}
		0=\int_X ud\mu=\int_X (u-c)d\mu+c\mu(X)=
		\int_\Omega u_+d\sigma+\int_\Omega -u_-d\sigma+c\mu(X)
	\end{equation*}
	Hence
	\begin{equation*}
		\int_X u_-d\mu=\int_\Omega u_+d\sigma+c\mu(X)\geq 3c\sigma(\Omega)\,.
	\end{equation*}
	and regarding $c-u_-$,
	\begin{equation*}
		\int_{\Omega}(c-u_-)d\sigma\leq -2c\sigma(\Omega)\,.
	\end{equation*}
	Since $\Psi$ is bilipschitz conformal, we can apply Lemma~\ref{lem:bettergradient}
	and we obtain that there is a constant $C(\Psi,\beta)$ such that
	\begin{equation*}
		\|\nabla u_-\|_m^m\geq C(\Psi,\beta)c^m\sigma(\Omega)\,.
	\end{equation*}
	The control on $\|\nabla u_+\|_m$ is then obtained from
	the Polya--Szegö inequality, asserting that
	\begin{equation*}
		\|\nabla u_+\|_m^m+\|\nabla u_-\|_m^m\leq 1\,.
	\end{equation*}
\end{proof}
\textbf{Step 3: controling $u_+$}
\begin{proposition}\label{prop:step3}
	With the notations introduced, there is a constant $C_{10}=C_{10}(V_0,\alpha,\varphi)$
	such that $u_+$ satisfies:
	\begin{equation}
		\int_\Omega F(c+u_+)d\sigma\leq C_{10}\,.
	\end{equation}
\end{proposition}
\begin{proof}
	Let $C_5$ denote the constant from Proposition~\ref{prop:step2}.
	Let $R$ be the constant:
	\begin{equation*}
		R=(1-C_6c^m\sigma(\Omega))^\frac{1}{m-1}
	\end{equation*}
	Wihtout loss of generality, we assume that $\frac{1}{2}\leq R<1$.
	Letting $v=\frac{u_+}{R^\frac{m-1}{m}}$, we get that $v$ is compactly supported
	and $\|\nabla v\|_m\leq 1$. Hence there is a constant $C_7(\varphi)$ such that
	\begin{equation*}
		\int_\Omega F(\alpha v^\frac{m}{m-1})d\sigma\leq C_7(\varphi)\,.
	\end{equation*}
	A quick computation shows that the function $h$
	\begin{equation*}
		h(t)=\frac{t^\frac{m}{m-1}}{R}-(t+c)^\frac{m}{m-1}
	\end{equation*}
	is lower bounded on $\R_+$, and its minimum is the value
	\begin{equation*}
		\min h=\frac{c^\frac{m}{m-1}}{(1-R^{m-1})^\frac{1}{m-1}}
			=\frac{1}{(C_6\sigma(\Omega))^\frac{1}{m-1}}\,.
	\end{equation*}
	In particular, we can bound the integral of $F(c+u_+)$:
	\begin{flalign*}
		\int_\Omega F(\alpha(c+u_+)^\frac{m}{m-1})d\sigma
			&\leq\int_\Omega \exp(\alpha(c+u_+)^\frac{m}{m-1})d\sigma\\
			&\leq e^{\frac{\alpha}{(C_3\sigma(\Omega))^\frac{1}{m-1}}}\int_\Omega
			\exp(\alpha v^\frac{m}{m-1})d\sigma\\
			&\leq  e^{\frac{\alpha}{(C_3\sigma(\Omega))^\frac{1}{m-1}}}\big(C_7
				+\sum_{j=0}^{m-2}\frac{\alpha^j}{j!}\int_\Omega |v|^\frac{jm}{m-1}d\sigma\big)
	\end{flalign*}
	Since $1\leq\frac{jm}{m-1}\leq m$, we can dominate $x^\frac{jm}{m-1}$ by $1+x^m$
	to obtain bounds $C_8(m,\alpha,\varphi)$ and $C_9(m,\alpha,\varphi)$
	such that
	\begin{flalign*}
			\int_\Omega F(\alpha(c+u_+)^\frac{m}{m-1})d\sigma\leq
			e^{\frac{\alpha}{(C_5\sigma(\Omega))^\frac{1}{m-1}}}
			(C_8+C_9\sigma(\Omega))\,.
	\end{flalign*}
	This right-hand side, seen as a function of $\sigma(\Omega)$, is upper bounded on
	$(V_0,+\infty)$. Hence there is a bound $C_{10}(m,V_0,\alpha,\varphi)>0$
	such that
	\begin{equation*}
		\int_\Omega F(c+u_+)d\sigma\leq C_{10}\,,
	\end{equation*}
	as wished.
\end{proof}
\begin{proof}[Proof of theorem~\ref{thm:MTcompact}]
With the notations introduced, we compute the integral:
	\begin{equation*}
		\int_X F(\alpha u^\frac{m}{m-1})d\mu
		=\int_\Omega F(c+u_+)d\sigma+\int_\Omega F(c-u_-)d\sigma
	\end{equation*}
	From Propositions~\ref{prop:step1} and \ref{prop:step3} we obtain that
	\begin{equation*}
		\int_X F(\alpha u^\frac{m}{m-1})d\mu\leq C_5(\varphi,\alpha,V_0,\lambda_m)+C_{10}(\alpha,\varphi,V_0)\,,
	\end{equation*}
	as claimed.
\end{proof}

\newpage
\bibliographystyle{alpha}
\bibliography{references}

\begin{thebibliography}{APPS24}

\bibitem[AL89]{AL89}
Frederick J.~jun. Almgren and Elliott~H. Lieb.
\newblock Symmetric decreasing rearrangement is sometimes continuous.
\newblock {\em Journal of the American Mathematical Society}, 2(4):683--773,
  1989.

\bibitem[APP22]{APP22}
Gioacchino Antonelli, Enrico Pasqualetto, and Marco Pozzetta.
\newblock Isoperimetric sets in spaces with lower bounds on the {Ricci}
  curvature.
\newblock {\em Nonlinear Analysis. Theory, Methods \& Applications. Series A:
  Theory and Methods}, 220:59, 2022.
\newblock Id/No 112839.

\bibitem[APPS24]{APP+24}
Gioacchino Antonelli, Enrico Pasqualetto, Marco Pozzetta, and Daniele Semola.
\newblock Asymptotic isoperimetry on non collapsed spaces with lower {Ricci}
  bounds.
\newblock {\em Mathematische Annalen}, 389(2):1677--1730, 2024.

\bibitem[BKT24]{BKT24}
Zolt{\'a}n~M. Balogh, Alexandru Krist{\'a}ly, and Francesca Tripaldi.
\newblock Sharp log-{S}obolev inequalities in {CD}(0, {N})) spaces with
  applications.
\newblock {\em Journal of Functional Analysis}, 286(2):41, 2024.
\newblock Id/No 110217.

\bibitem[BM82]{BM82}
Pierre B{\'e}rard and Daniel Meyer.
\newblock In{\'e}galit{\'e}s isop{\'e}rim{\'e}triques et applications.
  ({Isoperimetric} inequalities and applications).
\newblock {\em Annales Scientifiques de l'{\'E}cole Normale Sup{\'e}rieure.
  Quatri{\`e}me S{\'e}rie}, 15:513--541, 1982.

\bibitem[Bro23]{Bro23H4}
Samuel Bronstein.
\newblock Almost-fuchsian structures on disk bundles over a surface.
\newblock {\em arXiv preprint arXiv:2305.06665}, 2023.

\bibitem[Bus82]{Bus82}
Peter Buser.
\newblock A note on the isoperimetric constant.
\newblock {\em Ann. Sci. {\'E}c. Norm. Sup{\'e}r. (4)}, 15:213--230, 1982.

\bibitem[BZ88]{BZ88}
John~E. Brothers and William~P. Ziemer.
\newblock Minimal rearrangements of {Sobolev} functions.
\newblock {\em Journal f{\"u}r die Reine und Angewandte Mathematik},
  384:153--179, 1988.

\bibitem[DEJ14]{DEJ14}
Jean Dolbeault, Maria~J Esteban, and Gaspard Jankowiak.
\newblock Onofri inequalities and rigidity results.
\newblock {\em arXiv preprint arXiv:1404.7338}, 2014.

\bibitem[DPG18]{DG18}
Guido De~Philippis and Nicola Gigli.
\newblock Non-collapsed spaces with {Ricci} curvature bounded from below.
\newblock {\em Journal de l'{\'E}cole Polytechnique -- Math{\'e}matiques},
  5:613--650, 2018.

\bibitem[DPM21]{DM21}
Nicol{\`o} De~Ponti and Andrea Mondino.
\newblock Sharp {C}heeger-{B}user type inequalities in
  ({$\mathrm{RCD}(K,\infty)$}) spaces.
\newblock {\em The Journal of Geometric Analysis}, 31(3):2416--2438, 2021.

\bibitem[Gal88]{Gal88}
Sylvestre Gallot.
\newblock In{\'e}galit{\'e}s isop{\'e}rim{\'e}triques et analytiques sur les
  vari{\'e}t{\'e}s riemanniennes.
\newblock {\em Ast{\'e}risque}, 163(164):31--91, 1988.

\bibitem[Gig15]{Gig15}
Nicola Gigli.
\newblock {\em On the differential structure of metric measure spaces and
  applications}, volume 1113 of {\em Memoirs of the American Mathematical
  Society}.
\newblock Providence, RI: American Mathematical Society (AMS), 2015.

\bibitem[Haj96]{Haj96}
Piotr Haj{\l}asz.
\newblock Sobolev spaces on an arbitrary metric space.
\newblock {\em Potential Analysis}, 5(4):403--415, 1996.

\bibitem[Heb96]{Heb96}
Emmanuel Hebey.
\newblock {\em Sobolev spaces on Riemannian manifolds}, volume 1635.
\newblock Springer Science \& Business Media, 1996.

\bibitem[Kri19]{Kri19}
Alexandru Krist{\'a}ly.
\newblock New geometric aspects of {Moser}-{Trudinger} inequalities on
  {Riemannian} manifolds: the non-compact case.
\newblock {\em Journal of Functional Analysis}, 276(8):2359--2396, 2019.

\bibitem[Li01]{Li01}
Yuxiang Li.
\newblock Moser-trudinger inequality on compact riemannian manifolds of
  dimension two.
\newblock {\em Journal of Partial Differential Equations}, 14(2):163--192,
  2001.

\bibitem[Li05]{Li05}
Yuxiang Li.
\newblock Extremal functions for the moser-trudinger inequalities on compact
  riemannian manifolds.
\newblock {\em Science China Mathematics}, 5(48):618--648, 2005.

\bibitem[LL21]{LL21}
Jungang Li and Guozhen Lu.
\newblock Critical and subcritical {Trudinger}-{Moser} inequalities on complete
  noncompact {Riemannian} manifolds.
\newblock {\em Advances in Mathematics}, 389:36, 2021.
\newblock Id/No 107915.

\bibitem[LMS10]{LMS10}
Barnab{\'e}~Pessoa Lima, Jos{\'e} F{\'a}bio~Bezerra Montenegro, and
  Newton~Lu{\'i}s Santos.
\newblock Eigenvalue estimates for the {(p)}-{Laplace} operator on manifolds.
\newblock {\em Nonlinear Analysis. Theory, Methods \& Applications. Series A:
  Theory and Methods}, 72(2):771--781, 2010.

\bibitem[LV07]{LV07}
John Lott and C{\'e}dric Villani.
\newblock Weak curvature conditions and functional inequalities.
\newblock {\em Journal of Functional Analysis}, 245(1):311--333, 2007.

\bibitem[Mat00]{Mat00}
Ana-Maria Matei.
\newblock First eigenvalue for the {$p$}-{Laplace} operator.
\newblock {\em Nonlinear Analysis. Theory, Methods \& Applications. Series A:
  Theory and Methods}, 39(8):1051--1068, 2000.

\bibitem[Mir03]{Mir03}
Michele Miranda.
\newblock Functions of bounded variation on ``good'' metric spaces.
\newblock {\em Journal de Math{\'e}matiques Pures et Appliqu{\'e}es.
  Neuvi{\`e}me S{\'e}rie}, 82(8):975--1004, 2003.

\bibitem[Mos71]{Mos71}
J{\"u}rgen Moser.
\newblock A sharp form of an inequality by n. trudinger.
\newblock {\em Indiana University Mathematics Journal}, 20(11):1077--1092,
  1971.

\bibitem[MS10]{MS10}
G~Mancini and K~Sandeep.
\newblock Moser--trudinger inequality on conformal discs.
\newblock {\em Communications in Contemporary Mathematics}, 12(06):1055--1068,
  2010.

\bibitem[MS20]{MS20}
Andrea Mondino and Daniele Semola.
\newblock Polya-{Szego} inequality and {Dirichlet} p-spectral gap for
  non-smooth spaces with {Ricci} curvature bounded below.
\newblock {\em Journal de Math{\'e}matiques Pures et Appliqu{\'e}es.
  Neuvi{\`e}me S{\'e}rie}, 137:238--274, 2020.

\bibitem[MST13]{MST13}
Gianni Mancini, Kunnath Sandeep, and Cyril Tintarev.
\newblock Trudinger-{Moser} inequality in the hyperbolic space {$(\mathbb
  H^N)$}.
\newblock {\em Advances in Nonlinear Analysis}, 2(3):309--324, 2013.

\bibitem[NV24]{NV24}
Francesco Nobili and Ivan~Yuri Violo.
\newblock Fine {P}{\'o}lya-{S}zegö rearrangement inequalities in metric spaces
  and applications.
\newblock Preprint, {arXiv}:2409.14182 [math.{AP}] (2024), 2024.

\bibitem[Raj12]{Raj12}
Tapio Rajala.
\newblock Local {P}oincar{\'e} inequalities from stable curvature conditions on
  metric spaces.
\newblock {\em Calculus of Variations and Partial Differential Equations},
  44(3-4):477--494, 2012.

\bibitem[Sch43]{Sch43}
Erhard Ing~Grad Schmidt.
\newblock Beweis der isoperimetrischen eigenschaft der kugel im hyperbolischen
  und sph{\"a}rischen raum jeder dimensionenzahl.
\newblock {\em Mathematische Zeitschrift}, 49:1--109, 1943.

\bibitem[Stu06a]{Stu06I}
Karl-Theodor Sturm.
\newblock On the geometry of metric measure spaces. {I}.
\newblock {\em Acta Mathematica}, 196(1):65--131, 2006.

\bibitem[Stu06b]{Stu06II}
Karl-Theodor Sturm.
\newblock On the geometry of metric measure spaces. {II}.
\newblock {\em Acta Mathematica}, 196(1):133--177, 2006.

\bibitem[Tru67]{Tru67}
Neil~S. Trudinger.
\newblock On {I}mbeddings into {O}rlicz {S}paces and {S}ome {A}pplications.
\newblock {\em Journal of Mathematics and Mechanics}, 17(5):473--483, 1967.

\bibitem[Vil09]{Vil09}
C{\'e}dric Villani.
\newblock {\em Optimal transport. {Old} and new}, volume 338 of {\em
  Grundlehren der Mathematischen Wissenschaften}.
\newblock Berlin: Springer, 2009.

\bibitem[Yan12]{Yan12}
Yunyan Yang.
\newblock Trudinger-{Moser} inequalities on complete noncompact {Riemannian}
  manifolds.
\newblock {\em J. Funct. Anal.}, 263(7):1894--1938, 2012.

\end{thebibliography}
\end{document}